\documentclass[12pt]{amsart}
\usepackage{amscd}
\usepackage{amssymb}
\usepackage{graphicx}
\usepackage{color}
\usepackage[centering,text={15.5cm,22cm},
		marginparwidth=20mm]{geometry}
\usepackage{srcltx}
\usepackage[]{hyperref}

\usepackage{mathrsfs}

\usepackage[all]{xy}

\newtheorem{theorem}{Theorem}[section]
\newtheorem{lemma}[theorem]{Lemma}
\newtheorem{proposition}[theorem]{Proposition}

\newtheorem{conjecture}[theorem]{Conjecture}

\theoremstyle{definition}
\newtheorem{definition}[theorem]{Definition}

\newtheorem{example}[theorem]{Example}

\theoremstyle{remark}

\numberwithin{equation}{section}
\numberwithin{figure}{section}

\newcommand{\NN} {\mathbb{N}}
\newcommand{\ZZ} {\mathbb{Z}}
\newcommand{\VV} {\mathbb{V}}
\newcommand{\QQ} {\mathbb{Q}}
\newcommand{\RR} {\mathbb{R}}
\newcommand{\CC} {\mathbb{C}}

\newcommand{\PP} {\mathbb{P}}
\renewcommand{\AA} {\mathbb{A}}
\newcommand{\GG} {\mathbb{G}}

\newcommand {\shC}  {\mathcal{C}}

\newcommand {\shL}  {\mathcal{L}}

\newcommand {\shN}  {\mathcal{N}}

\newcommand {\foD}  {\mathfrak{D}}

\newcommand {\foL}  {\mathfrak{L}}
\newcommand {\foM}  {\mathfrak{M}}

\newcommand {\foX}  {\mathfrak{X}}

\newcommand {\fod}  {\mathfrak{d}}

\newcommand {\fog}  {\mathfrak{g}}
\newcommand {\foh}  {\mathfrak{h}}

\newcommand {\fom}  {\mathfrak{m}}

\newcommand {\Ancestors}  {\operatorname{Ancestors}}

\newcommand {\as}  {\mathrm{as}}

\newcommand {\Aut}  {\operatorname{Aut}}

\newcommand {\Child} {\operatorname{Child}}

\newcommand {\ev}  {\operatorname{ev}}

\newcommand {\GL}  {\operatorname{GL}}

\newcommand {\hol}  {\mathrm{hol}}
\newcommand {\Hom}  {\operatorname{Hom}}

\newcommand {\id}  {\operatorname{Id}}

\newcommand {\Init}  {\operatorname{Init}}

\newcommand {\kk} {\CC}

\newcommand {\Leaves}  {\operatorname{Leaves}}

\newcommand {\length} {\mathrm{length}}

\newcommand {\lra}  {\longrightarrow}

\newcommand {\M} {\mathcal{M}}

\newcommand {\Mult} {\operatorname{Mult}}

\newcommand {\out}  {\mathrm{out}}

\renewcommand{\P}  {\mathscr{P}}

\newcommand {\Parents} {\operatorname{Parents}}

\newcommand {\PGL}  {\operatorname{PGL}}

\newcommand {\Scatter} {\operatorname{{\mathsf{S}}}}

\newcommand {\Sing} {\operatorname{Sing}}

\newcommand {\Spec} {\operatorname{Spec}}

\newcommand {\trop} {\mathrm{trop}}

\def\mapright#1{\smash{
  \mathop{\longrightarrow}\limits^{#1}}}

\def\mapdown#1{\Big\downarrow
   \rlap{$\vcenter{\hbox{$\scriptstyle#1$}}$}}

\def\mydate{\ifcase\month \or January\or February\or March\or
April\or May\or June\or July\or August\or September\or October\or 
November\or December\fi \space\number\day,\space\number\year}

\begin{document}

\title[The Tropical Vertex]{The Tropical Vertex}

\author{M. Gross} \address{Department of Mathematics, UCSD,
La Jolla, CA 92093, USA}
\email{mgross@math.ucsd.edu}

\author{R. Pandharipande}\address{Department of Mathematics, Princeton 
Univerity, Princeton, NJ 08544, USA}
\email{rahulp@math.princeton.edu}

\author{B. Siebert} \address{Department Mathematik,
Universit\"at Hamburg, 
20146 Hamburg,
Germany}
\email{bernd.siebert@math.uni-hamburg.de}

\date{February 4th, 2009}

\begin{abstract}
Elements of the tropical vertex group are formal families
of symplectomorphisms of the 2-dimensional algebraic
torus. 
We prove ordered product factorizations in the tropical vertex group are
equivalent to calculations of certain genus 0
relative Gromov-Witten invariants of toric surfaces.
The relative invariants which arise have full tangency  
to a toric divisor at a single unspecified point.
The method uses scattering diagrams, tropical curve counts,
degeneration formulas, and exact multiple cover
calculations in orbifold Gromov-Witten theory.
\end{abstract}

\maketitle
\tableofcontents

\section*{Introduction}
\subsection{Automorphisms of the torus}
The 2-dimensional complex torus has very few
automorphisms 
$$\theta: \CC^* \times \CC^* \rightarrow \CC^* \times \CC^*$$
as an algebraic group. Since $\theta$ must take each component
$\CC^*$ to a
$1$-dimensional
subtorus, 
$$\text{Aut}_\CC^{\text{Gr}}(\CC^* \times \CC^*) 
\stackrel{\sim}{=} \GL_2(\mathbb{Z}).$$
As a complex algebraic variety, $\CC^* \times \CC^*$ has, in
addition, only the automorphisms obtained by the translation action
on itself,{\footnote{We leave the elementary proof
to the reader. An argument can be found by using
the characterization
$$\phi(z) = \lambda \cdot z^k\, \ \ \ \  \lambda\in \CC^*, \ k \in \mathbb{Z}$$
of all algebraic maps $\phi: \CC^* \rightarrow \CC^*$.
}} 
$$1 \rightarrow \CC^* \times \CC^*
 \rightarrow 
\text{Aut}_{\CC}(\CC^* \times \CC^*) \rightarrow
\text{Aut}_\CC^{\text{Gr}}(\CC^* \times \CC^*)
  \rightarrow 1.$$

A much richer algebraic structure appears if formal 1-parameter
families of automorphisms of $\CC^*\times \CC^*$ 
are considered,
$$A = \text{Aut}_{\CC[[t]]}(\CC^* \times \CC^* \times \text{Spec}(\CC[[t]])).$$
Let $x$ and $y$ be the coordinates of the two factors
of $\CC^* \times \CC^*$. Then,
$$\CC^* \times \CC^* = \text{Spec}(\CC[x,x^{-1},y,y^{-1}]).$$
We may alternatively view $A$ as
a group of algebra automorphisms,
$$A = \text{Aut}_{\CC[[t]]} (\CC[x,x^{-1},y,y^{-1}][[t]]).$$

Nontrivial elements of $A$ are easily found.
Let $(a,b)\in \mathbb{Z}^2$ be a nonzero vector, and let
$f\in \CC[x,x^{-1},y,y^{-1}][[t]]$ be a function of the form
$$f= 1 + t x^a y^b\cdot g( x^ay^b,t), \ \ \ \ g(z,t) \in \CC[z][[t]].$$
We specify the values of an automorphism on $x$ and $y$ by
\begin{equation}\label{vrzza}
\theta_{(a,b),f\ }(x) = f^{-b}\cdot x, \ 
\ \theta_{(a,b),f\ }(y)= f^a\cdot y\ \ .
\end{equation}
The assignment \eqref{vrzza} extends uniquely to determine
an element $\theta_{(a,b),f} \in A$ with
$$\theta^{-1}_{(a,b),f}=\theta_{(a,b),f^{-1}}.$$

\subsection{Tropical vertex group}
The tropical vertex group $H\subseteq A$ is the completion
with respect to the maximal ideal $(t) \subseteq \CC[[t]]$
of the subgroup 
generated by {\em all} elements of the form $\theta_{(a,b),f}$.
 A more natural characterization of $H$
via the associated Lie algebra is reviewed in Section  \ref{tvg}.

The torus $\CC^* \times \CC^*$ has a standard 
holomorphic symplectic form given by 
$$\omega = \frac{dx}{x} \wedge \frac{dy}{y}.$$
Let $S\subseteq A$ be the subgroup of automorphisms
preserving $\omega$,
$$S = \{ \ \theta \in A\ | \ \theta^*(\omega)= \omega \ \}.$$ 
A direct calculation shows $H\subseteq S$.

A slight variant of the tropical vertex group
 $H$ first arose in the study of affine
structures by Kontsevich and Soibelman in \cite{ks}.
Further development, related to mirror symmetry and tropical geometry, can be found in 
\cite{GS}. 
Recently, the tropical vertex group has played a role in 
wall-crossing formulas for counting invariants in derived
categories \cite{ks2}.

\subsection{Commutators}
The first question  we can ask about the tropical
vertex group is to find
a formula for the commutators of the generators.
The answer, the main result of the paper,
turns out to be surprisingly subtle. The commutator formula
 is
expressed in terms of the relative Gromov-Witten
theories of toric surfaces. 

Perhaps the simplest nontrivial case to consider is the commutator of 
the elements
$$S_{\ell_1}=\theta_{(1,0),(1+tx)^{\ell_1}} \ \  \text{and} \ \ 
T_{\ell_2}=\theta_{(0,1), (1+ty)^{\ell_2}}$$
where $\ell_1,\ell_2 >0$.
By an elementary result of \cite{ks} reviewed in Section \ref{pop}, there exists a unique
factorization
\begin{equation}\label{jjtr}
    T_{\ell_2}^{-1} \circ S_{\ell_1} \circ T_{\ell_2}  \circ S_{\ell_1}^{-1} =
\stackrel{\rightarrow} \prod \theta_{(a,b),f_{(a,b)}}\ 
\end{equation}
where the product on the right is over {\em all} primitive
vectors $(a,b)\in \mathbb{Z}^2$ lying strictly in the first quadrant.{\footnote{A vector
 $(a,b)$ is primitive if it is not divisible in $\mathbb{Z}^2$. 
Primitivity implies $(a,b)\neq (0,0)$.
Strict inclusion in
the first  quadrant is equivalent to $a>0$ and $b>0$.}}
The order is determined by 
increasing slopes of the vectors $(a,b)$.
%
The question is what are the functions $f_{(a,b)}$ 
associated to the slopes?

\subsection{Toric surfaces} \label{ttss}

Let $(a,b)\in \mathbb{Z}^2$ be a primitive vector lying 
strictly in the first quadrant.
The rays generated by $(-1,0)$, $(0,-1)$, and $(a,b)$ determine
a complete rational fan in $\mathbb{R}^2$. Let $X_{a,b}$
be the associated toric surface{\footnote{$X_{a,b}$
is simply a weighted projective plane.
Arbitrary toric surfaces
will arise in the study of more general commutators.
}} with toric divisors
$$D_1, D_2, D_\out \subseteq X_{a,b}$$
corresponding to the respective rays.
Let 
$$ X_{a,b}^o \subseteq X_{a,b}$$
be the open surface obtained by removing the three
toric fixed points. Let 
$D_1^o, D_2^o, D_\out^o$ be the restrictions of the
toric divisors to $X_{a,b}^o$.

We denote {\em ordered 
partitions} ${\bf Q}$ of length $\ell$ by  $q_1+ \ldots + q_{\ell}$. 
Ordered partitions differ from usual partitions in two basic ways.
First,
the ordering of the parts matters.
Second, the parts $q_i$ are
required only to be non-negative integers (0 is permitted).
The {size}  $|{\bf Q}|$ is the sum of the parts.

Let $k\geq 1$.
Let 
${\bf P}_a=p_1+ \ldots+ p_{\ell_a}$ and ${\bf P}_b=p'_1 + \ldots +p'_{\ell_b}$ 
be ordered partitions of size $a k$ 
and $bk$ respectively. Denote the pair by ${\bf{P}}=({\bf P}_a,{\bf P}_b)$.
Let
$$\nu: X_{a,b}[{\bf P}] \rightarrow X_{a,b}$$
be the blow-up of $X_{a,b}$ along $\ell_a$ and $\ell_b$
distinct points of $D^o_1$ and $D^o_2$.
Let 
$$X^o_{a,b}[{\bf P}] = \nu^{-1}(X^o_{a,b}).$$

Let $\beta_k \in H_2(X_{a,b},\mathbb{Z})$ be the unique
class with intersection numbers
$$\beta_k \cdot D_1 =ak , \ \ \beta_k \cdot D_2 = bk, \ \ 
\beta_k\cdot D_\out =k.$$
Let $E_i$ and $E'_j$ be the $i^{th}$ and $j^{th}$
exceptional divisors over $D^o_1$ and $D^o_2$.
Let $$\beta_k[{\bf P}] = \nu^*(\beta_k) -\sum_{i=1}^{\ell_a} p_i [E_i]
-\sum_{j=1}^{\ell_b} p'_j [E'_j]\ \in H_2(X_{a,b}[{\bf P}], \mathbb{Z}) .$$

By a parameter count, the moduli space 
$\overline{\foM}(X^o_{a,b}[{\bf P}]/D_\out^o)$
of genus 0 
maps to  $X^o[{\bf P}]$ in class $\beta_k[{\bf P}]$
with full contact order  $k$ at an unspecified point of $D^o_\out$
is of virtual dimension 0.
 In Section \ref{gwt}, we will show the 
corresponding Gromov-Witten invariant
$$N_{a,b}[{\bf P}] \in \mathbb{Q}$$
is well-defined.{\footnote{The only issue  
to understand is the relative Gromov-Witten theory of the
{\em open} geometry $X_{a,b}^o[{\bf P}]/D_\out^o$. 
We will
show there is sufficient properness here to define Gromov-Witten
invariants in the usual way.}}

\subsection{Formula}
Since the series $f_{(a,b)}$ in \eqref{jjtr}
 starts with 1, we may take the logarithm.
Homogeneity constraints determine the behavior of the
variable $t$. We define the coefficients
$c^k_{a,b}(\ell_1,\ell_2) \in \mathbb{Q}$ by
$$\log f_{(a,b)} = \sum_{k\geq 1}  k\ c^k_{a,b}(\ell_1,\ell_2) \cdot  
(tx)^{ak}\ (ty)^{bk} .$$
The function $f_{(a,b)}$ is linked to Gromov-Witten theory
by the following result.

\begin{theorem}\label{fff}
We have
$$c^k_{a,b}(\ell_1,\ell_2) = \sum_{|{\bf P}_a| = ak} \ \sum_{|{\bf P}_b|=bk} 
N_{a,b}[({\bf P}_a,{\bf P}_b)]$$
where the sums are over all ordered
partitions ${\bf P}_a$ of size {\em ak} and length $\ell_1$
  and ${\bf P}_b$
of size {\em bk} and length $\ell_2$. 
\end{theorem}

Theorem \ref{fff} is the first in a series of results. 
A refinement of  Theorem \ref{fff} proven
in Section \ref{fcf} shows the invariants $N_{a,b}[{\bf P}]$
are {\em determined} by ordered product factorizations of
commutators in the tropical vertex group.
In fact, the tropical vertex group calculations are much 
 simpler (and much more conceptually appealing)
than the alternative methods available
for calculating $N_{a,b}[{\bf P}]$.

Commutator formulas for arbitrary generators
$\theta_{(a,b),f}$ and $\theta_{(a',b'),f'}$
of the tropical vertex group are proven in Section \ref{ftvg}.
For general functions $f$ and $f'$, orbifold blow-ups of
toric surfaces replace the ordinary blow-ups discussed
above. As a result, the relations between generators of
the tropical vertex group are completely described in terms
of Gromov-Witten invariants.

The tropical vertex group was used in  \cite{ks} to
construct
 rigid analytic K3 surfaces from affine 2-spheres with singularities.
The group was similarly used in \cite{GS} to construct 
explicit degenerations of Calabi-Yau manifolds from integral affine
manifolds with singularities. The latter work gives a precise description
of the B-model for Calabi-Yau manifolds (as well as manifolds with
effective anti-canonical bundle).
The commutator formula of   
Theorem \ref{fff} provides an interpretation of 
the tropical vertex group in terms of curve counts on the
A-model side.
Hence, Theorem 0.1 should be viewed as a mirror
symmetry relation. 

Special cases of commutators in the tropical vertex group
have direct interpretation in the wall-crossing work of \cite{ks2,Re}.
The $$\ell=\ell_1=\ell_2$$ case of Theorem 0.1 is related to the
quiver $Q_\ell$ with two vertices and  $\ell$ 
arrows (in the same direction).
Reineke \cite{Re} has proven  the functions $f_{(a,b)}$ 
are then determined by the Euler characteristics of the moduli spaces
of semistable representations of $Q_\ell$ with dimension
vectors along the ray generated by $(a,b)$. Whether such
alternative interpretations hold for more general
commutators is an interesting question.

\subsection{Plan of the paper}
We start with a formal discussion of the tropical vertex group and
 the associated scattering diagrams for path oriented products
in Section \ref{sctt}. The scattering diagram expansions connect
the commutators to tropical curve counts explained in Section \ref{tcc}.
The tropical counts are related to the enumeration of holomorphic
curves in Section 
\ref{section-tropicaltoholomorphic}
 following \cite{Mk,nisi}. The Gromov-Witten theory
of open toric surfaces is discussed in Section \ref{gwt}.
The commutator formulas are proven in Section \ref{ftvg} via
degeneration and exact Gromov-Witten calculations.

The Gromov-Witten invariants $N_{a,b}[{\bf P}]$ which arise in Theorem \ref{fff} are virtual
counts with complicated multiple cover contributions of excess dimension.
The corresponding (conjectural) BPS structure  
is discussed at the end of the
paper in Section \ref{bps}. The method also applies to neighboring questions
such as the enumeration of degree $d$ rational curves
in $\mathbb{P}^2$ with full tangency at a single point 
to a nonsingular elliptic curve $E\subseteq \mathbb{P}^2$ studied
in \cite{Ga,Ta}.

\subsection{Acknowledgments}
We thank J. Bryan, T. Graber, J. de Jong, J. Koll\'ar, M.
Kontsevich, D. Maulik, Y. Soibelman, and 
N. Takahashi for many related comments and discussions.
M.~G. was partially supported by NSF grant DMS-0805328 and
the DFG.
R.~P. was partially supported by NSF grant DMS-0500187 and the
Clay foundation.

\section{Scattering} \label{sctt}
\subsection{Tropical vertex group (again)}\label{tvg}
We begin by reviewing ideas from \cite{GS}. While most of
the material we need can be extracted from \cite{GS},
the development there is in much greater generality  than we need,
and a number of points simplify in our context. 
We will give a completely self-contained exposition.

We fix once and for all a lattice $M\cong \ZZ^2$ with basis $e_1=(1,0)$ and
$e_2=(0,1)$. Let
  $$N=\Hom_{\ZZ}(M,\ZZ),\ \ 
M_{\RR}=M\otimes_{\ZZ}\RR,\ \ N_{\RR}=N\otimes_{\ZZ}\RR \ .$$ 
%
For $m\in M$, let $z^m\in \kk[M]$ denote the
corresponding element in the group ring. Let
$$x= z^{e_1}, \ \ y=z^{e_2}.$$
Then,
$\kk[M]$ is simply  the ring of Laurent polynomials
in $x$ and $y$.

In what follows, let $R$ be an Artin local $\kk$-algebra or a complete
local $\kk$-algebra.{\footnote{In the introduction, we took $R=\CC[[t]]$.}}
Let $\fom_{R}\subseteq R$ be the maximal ideal. We define the
\emph{module of log derivations} of 
\[
\kk[M]\widehat\otimes_{\kk} R=
\lim_{\longleftarrow} \kk[M]\otimes_{\kk} R/\fom_R^k
\]
to be
\[
\Theta(\kk[M]\widehat\otimes_{\kk} R)=\Hom(M,\kk[M]\widehat\otimes_{\kk} R)=
(\kk[M]\widehat\otimes_{\kk} R)\otimes_{\ZZ} N.
\]
A log derivation $\xi$ induces an ordinary derivation $\bar{\xi}$ of 
$\kk[M]\widehat\otimes_{\kk} R$ over $R$ via the formula
\[
\bar\xi(z^m)=\xi(m)z^m.
\]
We will write $a\partial_n$ for 
$a\otimes n \in (\kk[M]\widehat\otimes_{\kk} R)\otimes_{\ZZ} N =
\Theta(\kk[M]\widehat\otimes_{\kk} R)$
and for the associated ordinary derivation,
\[
(a\partial_n)(z^m)=a\langle m,n\rangle z^m.
\]

Let $\fog_R=\fom_R\Theta(\kk[M]\widehat\otimes_{\kk} R)$.
Given any $\xi\in \fog_R$, we obtain
an element
\[
\exp(\xi)\in \Aut_R(\kk[M]\widehat\otimes_{\kk} R)
\]
of the group of ring automorphisms defined by
\[
\exp(\xi)(a)=\text{Id}(a)+\sum_{i=1}^{\infty} {\bar\xi^i(a)\over i!}.
\]
The series makes sense by the assumptions on $R$.

In fact, $\fog_R$ is a Lie algebra under the natural bracket
defined by
\begin{eqnarray}
\label{Liebracket}
[z^m\partial_n,z^{m'}\partial_{n'}]
&=&(z^m\partial_n(z^{m'}))\partial_{n'}-(z^{m'}\partial_{n'}(z^m))\partial_n
\nonumber\\
&=&z^{m+m'}(\langle m',n\rangle\partial_{n'}-\langle m,n'\rangle
\partial_n)\\
&=&z^{m+m'}\partial_{\langle m',n\rangle n'-\langle m,n'\rangle n}.
\nonumber
\end{eqnarray}
By the Baker-Campbell-Hausdorff formula, the subset
$$G_R=\{\exp(\xi)\,|\, \xi\in\fog_R\}$$
of $\Aut_{R}(\kk[M]\widehat\otimes_{\kk} R)$
is a subgroup.

The subspace $\foh_R\subseteq \fog_R$ defined by
\[
\foh_R=\bigoplus_{m\in M\setminus\{0\}}
z^m( \fom_R \otimes m^{\perp})
\]
is closed under the Lie bracket, and hence via exponentiation
defines a subgroup
\[
\VV_R\subseteq G_R
\]
which we call the {\em tropical vertex group}.
The subgroup $\VV_R$ is  closed in the $\fom_R$-adic topology.

\subsection{Scattering diagrams}
We will depict certain products of automorphisms
of the form $\theta_{(a,b),f}$ as introduced in the Introduction
by diagrams in the plane $M_\mathbb{R}$. Automorphisms of the form
$\theta_{(a,b),f}$ in fact generate $\VV_R$, and this will give us a way
of describing elements in $\VV_R$.

\begin{definition}
 A \emph{ray} or \emph{line} is a pair $(\fod,f_{\fod})$ such that 
\begin{itemize}
\item $\fod\subseteq M_{\RR}$
is given by
\[
\fod=m_0'+\RR_{\ge 0} m_0
\]
if $\fod$ is a ray and
\[
\fod=m_0'+\RR m_0
\]
if $\fod$ is a line,
for some $m_0'\in M_{\RR}$ and $m_0\in M\setminus \{0\}$. The set
$\fod$ is called the \emph{support} of the line or ray. If $\fod$ is 
a ray, $m_0'$ is called the \emph{initial point} of the ray, written
as $\Init(\fod)$.
\item $f_{\fod}\in \kk[z^{m_0}]\widehat\otimes_{\kk} 
R\subseteq \kk[M]\widehat\otimes_{\kk} R.$
\item $f_{\fod}\equiv 1 \mod z^{m_0}\fom_R$.
\end{itemize}
\end{definition}

\begin{definition}
 A \emph{scattering diagram} $\foD$ is a set of lines and rays such
that for every power $k>0$, 
there are only a finite number of $(\fod,f_{\fod})\in \foD$
with $f_{\fod}\not\equiv 1 \mod \fom_R^k$.
\end{definition}

If $\foD$ is a scattering diagram, we write 
\[
\Sing(\foD)=\bigcup_{\fod\in\foD} \partial\fod \cup \bigcup_{\fod_1,\fod_2
\atop\dim \fod_1\cap\fod_2=0} \fod_1\cap \fod_2.
\]
Here $\partial\fod=\{\Init(\fod)\}$ if $\fod$ is a ray, and is empty
if $\fod$ is a line.

We are using a different sign
convention than in \cite{GS}. Here, we view a monomial $z^m$ as
propagating in the direction $m$, while in \cite{GS}, monomials propagate
in the direction $-m$. The negative sign is necessary in \cite{GS}. We
do not need the sign for our purposes here, so we dispense with it.

\subsection{Path ordered products}
\label{pop}
Consider a smooth immersion
$$\gamma:[0,1]\rightarrow  M_{\RR}\setminus\Sing(\foD)$$ 
with endpoints not contained in the scattering diagram $\foD$.
If all intersections of $\gamma$ with the rays and lines of  
$\foD$ are transverse,
 we can define $\theta_{\gamma,\foD}
\in \VV_R$, the {\em $\gamma$-ordered product of $\foD$} as follows.

For each power $k>0$, we can find
numbers 
\[
0<t_1\le t_2\le \cdots \le t_s<1
\]
and elements $\fod_i\in\foD$ with $f_{\fod_i}\not\equiv 1\mod \fom_{R}^k$
such that 
$$\gamma(t_i)\in\fod_i,$$
$\fod_i\not=\fod_j$ if $t_i=t_j$ and $i\not=j$, and $s$ taken as large as
possible. For each $i$, define $\theta_{\fod_i}$ to be 
\[
\theta_{\fod_i}=\exp(\log(f_{\fod_i})\partial_{n_0})
\]
where $n_0\in N$ is  primitive, annihilates the tangent
space to $\fod_i$, and is uniquely determined by the sign
convention 
\[
\langle n_0,\gamma'(t_i)\rangle >0.
\]
We define
\[
\theta^k_{\gamma,\foD}=\theta_{\fod_s}\circ\cdots\circ\theta_{\fod_1}.
\]
If $t_i=t_{i+1}$, then $\gamma(t_i)\in\fod_i\cap
\fod_{i+1}$. Hence, $\dim\fod_i\cap\fod_{i+1}=1$. The elements
$\theta_{\fod_i}$ and $\theta_{\fod_{i+1}}$ are easily 
seen to commute.
Therefore,  the product does not depend on any choices.
Let
$$\theta_{\gamma,\foD} = \lim_{k\rightarrow \infty}\theta^k_{\gamma,
\foD}\ .$$

\begin{example} 
\label{commutatoreg}
A basic example for us is given by the scattering diagram
\[
\foD=\{(\fod_1,f_1),(\fod_2,f_2)\},
\]
where $\fod_1$, $\fod_2$ are transverse lines through the origin.
A suitably chosen loop $\gamma$ around the origin produces a commutator:
\[
\theta_{\gamma,\foD}=\theta_2^{-1}\theta_1^{-1}\theta_2\theta_1\ ,
\]
where the elements 
$\theta_1$ and $\theta_2$ are obtained
 from $\fod_1$ and $\fod_2$ respectively,
see Figure \ref{commutatorfig}.
\begin{figure}
\input{commutator.pstex_t}
\caption{}
\label{commutatorfig}
\end{figure}
\end{example}

The following result was obtained by Kontsevich and Soibelman in \cite{ks}
with a number of variants proved in \cite{GS}. Since the proof is so simple,
we reproduce it here.

\begin{theorem}
\label{KSLemma}
Let $\foD$ be a scattering diagram. Then, there exists a scattering
diagram $\Scatter(\foD)$ containing $\foD$ such 
that $\Scatter(\foD)\setminus\foD$ consists
only of rays, and such that $\theta_{\gamma,\Scatter(\foD)}=\text{\em Id}$ 
for any closed
loop $\gamma$ for which $\theta_{\gamma,\Scatter(\foD)}$ is defined.
\end{theorem}

\proof We proceed inductively on $k$, showing that there exists a $\foD_k$
such that 
$$\theta_{\gamma,\foD_k}\equiv \id \mod \fom_R^{k+1}$$ for all
closed loops $\gamma$ for which $\theta_{\gamma,\foD_k}$ is defined.
We take $\foD_0=\foD$. To obtain $\foD_{k}$ from $\foD_{k-1}$, we proceed
as follows. Let $\foD'_{k-1}$ consist of those rays and lines $\fod$ in 
$\foD_{k-1}$ with $f_{\fod}\not\equiv 1\mod \fom_R^{k+1}$. This is a finite
set, so $\Sing(\foD'_{k-1})$ is finite also. 
Let $p\in\Sing(\foD'_{k-1})$. Let $\gamma_p$ be a closed simple 
loop around $p$, small enough so it contains no other points of 
$\Sing(\foD'_{k-1})$. Certainly, $$\theta_{\gamma_p,\foD_{k-1}}\equiv
\theta_{\gamma_p,\foD'_{k-1}}\mod \fom^{k+1}_R.$$ By the inductive
assumption we can write uniquely
\[
\theta_{\gamma_p,\foD'_{k-1}}=\exp(\sum_{i=1}^s c_iz^{m_i}\partial_{n_i})
\]
with $m_i\in M\setminus \{0\}$, $n_i\in m_i^{\perp}$ primitive and $c_i\in
\fom_R^k$. Let 
\[
\foD[p]=\{(p+\RR_{\ge 0}m_i,1\pm c_iz^{m_i})\,|\,i=1,\ldots,s\}.
\]
The sign is chosen in each ray so that its contribution 
to $\theta_{\gamma_p,\foD[p]}$ is $\exp(-c_iz^{m_i}\partial_{n_i})$ 
modulo $\fom_R^{k+1}$.
Since $c_iz^{m_i}\partial_{n_i}$ is in the center of
$\foh_{R/\fom_R^{k+1}}$, we see 
\[
\theta_{\gamma_p,\foD_{k-1}\cup\foD[p]}=\id \mod \fom_R^{k+1}.
\]
Similarly, any automorphism coming from $\foD[p]$ commutes with any
automorphism coming from $\foD_{k-1}$ modulo $\fom_R^{k+1}$, and hence
\[
\foD_k=\foD_{k-1}\cup \bigcup_{p} \foD[p]
\]
has the desired properties.

We take $\Scatter(\foD)$ to be the (non-disjoint) union of the
$\foD_k$'s. The diagram $\Scatter(\foD)$ will usually have
infinitely many rays even if $\foD$ has finitely many.
\qed
\vspace{10pt}

The scattering diagram $\Scatter(\foD)$ is not unique.
There are always trivial changes which can be made to scattering diagrams.

\begin{definition} Two scattering diagrams $\foD$, $\foD'$ over a ring
$R$ are \emph{equivalent} if $$\theta_{\gamma,\foD}=\theta_{\gamma,\foD'}$$
for every curve $\gamma$ for which both sides are defined.
\end{definition}

Every scattering diagram $\foD$ is equivalent to a unique
{\em minimal} scattering diagram which 
does not have two rays or lines
with the same support and contains no trivial rays or lines --- rays or lines
$(\fod,f_{\fod})$ with $f_{\fod}=1$. In fact,
if $(\fod_1,f_1),(\fod_2,f_2)\in\foD$
with $\fod_1=\fod_2$, then we can replace these two rays or lines
with the single
ray or line $(\fod_1,f_1f_2)$ without affecting any automorphism 
$\theta_{\gamma,\foD}$. In addition we can remove any trivial ray or line.
The diagram  $\Scatter(\foD)$ is easily seen to be unique up
to equivalence. We will often assume $\Scatter(\foD)$ is minimal.

The operation $\Scatter$ is functorial in the following sense. 
If $\varphi:R\rightarrow
R'$ is a ring homomorphism, and 
\[
\varphi(\foD)=\{(\fod,\varphi(f_{\fod}))\,|\, (\fod,f_{\fod})\in \foD\},
\]
then $\varphi(\Scatter(\foD))$ is equivalent to $\Scatter(\varphi(\foD))$.

\begin{example} 
\label{basicexample}
Let $R=\kk [[ t_1,t_2]]$ and 
consider for $\ell_1,\ell_2>0$ the diagram 
\[
\foD=\{(\RR (1,0),(1+t_1x)^{\ell_1}), (\RR (0,1),(1+t_2y)^{\ell_2})\}.
\]
Here, $x=z^{e_1}$ and $y=z^{e_2}$ as before.
Theorem \ref{KSLemma} produces a scattering diagram
$\Scatter(\foD)$. We will later see that all elements of
$\Scatter(\foD)\setminus\foD$ 
lie in the first quadrant. However, 
$\Scatter(\foD)\setminus\foD$ can be very complicated.{\footnote{The behavior
of these scattering diagrams has also been studied by Kontsevich.}} 

If $\ell_1=\ell_2=1$, then $\Scatter(\foD)$ requires only a single
additional ray,
\[
\Scatter(\foD)\setminus\foD=\{(\RR(1,1), 1+t_1t_2xy)\},
\]
as can be easily checked by hand. See Figure \ref{d1=1,d2=1}.
\begin{figure}
\input{oneone.pstex_t}
\caption{$\Scatter(\foD)$ for $\ell_1=1$, $\ell_2=1$. Here the automorphisms 
are given
explicitly, and the identity $\theta_{\gamma,\Scatter(\foD)}$ is just the composition of the
given automorphisms.}
\label{d1=1,d2=1}
\end{figure}

If $\ell_1=\ell_2=2$, then the following is found:
\begin{eqnarray*}
\Scatter(\foD)
\setminus\foD=&&\{(\RR (n+1,n),(1+(t_1x)^{n+1}(t_2y)^n)^2)| n\in\ZZ, n\ge 1\}\\
&\cup&
\{(\RR (n,n+1),(1+(t_1x)^n(t_2y)^{n+1})^2)| n\in\ZZ, n\ge 1\}\\
&\cup&\{(\RR (1,1), (1-t_1t_2xy)^{-4})\}.
\end{eqnarray*}

If $\ell_1=\ell_2=3$, then computer experiments (the algorithm
given in the proof of Theorem \ref{KSLemma} is easily
implemented) suggest the following behavior
for $\Scatter(\foD)$. First, as noticed by Kontsevich,
$\Scatter(\foD)$ has a certain periodicity. Namely,
\[
(\RR_{\ge 0} (m_1,m_2),f(z^{(m_1,m_2)}))
\in\Scatter(\foD)
\]
if and only if
\[
(\RR_{\ge 0} (3m_1-m_2,m_1),f(z^{(3m_1-m_2,m_1)}))
\in\Scatter(\foD),
\]
provided that $m_1,m_2$ and $3m_1-m_2$ are all positive. 
In addition, there are rays with support
$\RR_{\ge 0} (3,1)$ and $\RR_{\ge 0} (1,3)$, hence by the
 periodicity, there are also
rays with support 
$$\RR_{\ge 0}(8,3),\ \RR_{\ge 0}(21,8),\ \ldots \ \ \
\text{and}\  \ \  \RR_{\ge 0}(3,8),\ \RR_{\ge 0}(8,21),\ \ldots $$ 
which converge to
the rays of slope $(3\pm\sqrt{5})/2$, corresponding to the two distinct
eigenspaces of the linear transformation 
$\begin{pmatrix}3&-1\\1&0\end{pmatrix}$.
On the other hand, inside the cone generated by the rays of slope
$(3\pm\sqrt{5})/2$, it appears that 
every rational slope occurs. The function attached to the line of slope $1$
appears to be 
\[
\left(\sum_{n=0}^{\infty} {1\over 3n+1}\begin{pmatrix}4n\\ n\end{pmatrix}
(t_1x)^n(t_2y)^n\right)^9.
\]
\end{example}

\subsection{Asymptotic diagrams and factorization}
\label{asymfac}
A scattering diagram $\foD$ viewed from a great distance
yields an asymptotic diagram.

\begin{definition} The \emph{asymptotic
scattering diagram} $\foD_{\as}$ is obtained from
$\foD$ by replacing each ray $(m_0'+\RR_{\ge 0} m_0,f)$ with the ray
$(\RR_{\ge 0} m_0,f)$ emanating from the origin, and replacing each
line $(m_0'+\RR m_0,f)$ with the line $(\RR m_0,f)$ passing through
the origin. 
\end{definition}

If $\gamma$ is a sufficiently large simple
loop around the origin containing all points of $\Sing(\foD)$, then 
$$\theta_{\gamma,\foD}=\theta_{\gamma,\foD_{as}}.$$

We can often understand scattering diagrams
using a deformation technique as follows. 
Suppose $\foD$ consists of a number of lines $(\fod_i,f_i)$ through the origin,
and we wish to understand $\Scatter(\foD)$.
Assume that each $f_i$ can be factored as $f_i=\prod f_{ij}$. 
We replace
each line $(\fod_i,f_i)$ with the collection of lines 
$\{(m_{ij}'+\fod_i,f_{ij})\}$, where $m_{ij}'\in M_{\RR}$ is chosen generally.
Thus each line is split up into a number of parallel lines. Calling this
new scattering diagram $\widetilde\foD$, we can now apply Theorem \ref{KSLemma}
to obtain a scattering diagram $\Scatter(\widetilde\foD)$. 
Then $\Scatter(\widetilde\foD)_{\as}$
satisfies the conclusion of Theorem \ref{KSLemma} when applied to
$\foD$, since for a large simple loop around the origin,
$\theta_{\gamma,\Scatter(\widetilde\foD)_{as}}=\id$.

For example, in Example \ref{basicexample}, we can split the two
lines in $\foD$ into $\ell_1$ and $\ell_2$ parallel lines respectively.

\begin{definition} If $m\in M\setminus \{0\}$, 
the \emph{index} of $m$ is a positive integer $w$ such that $m=wm'$ for
$m'\in M$ primitive.
\end{definition}

The following simple calculation  is closely related to the notion of multiplicity
in tropical geometry.

\begin{lemma}
\label{nilpotentcollision}
Let $R=\kk[t_1,t_2]/(t_1^2,t_2^2)$, and let
\[
\foD=\{(\RR m_1,1+c_1t_1z^{w_1m_1}),(\RR m_2,1+c_2t_2z^{w_2m_2})\}
\]
where $m_1,m_2\in M$ are primitive with $m_1\wedge m_2\not=0$, 
$w_1,w_2$ are positive integers and
$c_1,c_2\in\kk$.
The scattering diagram $\Scatter(\foD)$ of 
Theorem \ref{KSLemma}
is obtained by adding the single ray
\begin{equation} \label{byy7}
\big(\RR_{\ge 0}(w_1m_1+w_2m_2),
1+c_1c_2t_1t_2w_{\out}|m_1\wedge m_2|z^{w_1m_1+w_2m_2}\big)
\end{equation}
where $|m_1\wedge m_2|$ denotes the absolute value of $m_1\wedge m_2\in
\bigwedge^2 M\cong\ZZ$, and $w_{\out}$ is the index of $w_1m_1+w_2m_2$.
\end{lemma}

If $w_1m_1+w_2m_2=0$  in \eqref{byy7}, no ray is added.
The associated function then is just 1 since then
$m_1\wedge m_2 =0$.

\proof
Labelling the two rays $\fod_1$ and $\fod_2$ respectively, choose a loop
$\gamma$ as depicted in Figure \ref{commutatorfig}, so that 
\[
\theta_{\gamma,\foD}=\theta_2^{-1}\theta_1^{-1}\theta_2\theta_1,
\]
where
\[
\theta_i=\exp(c_it_iz^{w_im_i}\partial_{n_i})
\]
with $n_i\in N$ primitive, orthogonal to $m_i$, $n_1$ negative on
$m_2$ and $n_2$ positive on $m_1$. The commutator is easily seen from
\eqref{Liebracket} or direct computation to be
\begin{equation}
\label{commexp}
\exp(c_1c_2t_1t_2z^{w_1m_1+w_2m_2}\partial_{w_1n_1\langle m_1,n_2\rangle
-w_2n_2\langle m_2,n_1\rangle}).
\end{equation}
Thus $\Scatter(\foD)$ needs to include exactly one additional ray, 
$(\RR_{\ge 0}(w_1m_1+w_2m_2),f)$, with $f$ chosen so
the corresponding 
automorphism determined by the crossing of
 $\gamma$ is the inverse
of \eqref{commexp}. 

First assume the basis $m_1,m_2$ of $M_{\RR}$
is positively oriented, as in Figure \ref{commutatorfig}. 
If we write
$m_i=(m_{i1},m_{i2})$ in some fixed basis of $M$, we can take 
$n_i=(m_{i2},-m_{i1})$ in the dual basis. Then,
$$\langle m_1,n_2\rangle=-\langle m_2,n_1\rangle=m_{11}m_{22}-m_{12}m_{21}
=|m_1\wedge m_2|$$
 because of the positive orientation. So
\begin{eqnarray*}
w_1n_1\langle m_1,n_2\rangle-w_2n_2\langle m_2,n_1\rangle
&=&
|m_1\wedge m_2|(w_1m_{12}+w_2m_{22},-w_1m_{11}-w_2m_{21})\\
&=&|m_1\wedge m_2|w_{\out}n_{\out}
\end{eqnarray*}
where $n_{\out}$ is a primitive element of $N$ pointing in the
same direction as 
$$w_1n_1\langle m_1,n_2\rangle-w_2n_2\langle m_2,n_1\rangle.$$
The automorphism we wish to attach to the new ray is
\[\exp(c_1c_2t_1t_2w_{\out}|m_1\wedge m_2|z^{w_1m_1+w_2m_2}
\partial_{-n_{\out}})\]
Noting that $-n_{\out}$ is positive on $m_2$, we see that we
can take
\[
f=1+c_1c_2t_1t_2w_{\out}|m_1\wedge m_2|z^{w_1m_1+w_2m_2}.
\] 
as claimed.

If $m_1,m_2$ is negatively oriented, the
argument is similar. Or, we can reduce
to the positively oriented case by switching labels.
\qed

\vspace{10pt}

In fact, using factorization and deformation,
Lemma \ref{nilpotentcollision}
 is enough to understand a scattering diagram
to any order. 
We will demonstrate the method for the scattering diagrams we 
require for the remainder of the paper, so the notation introduced here
will be used throughout.

\begin{definition} \label{ssss}

A {\em standard} scattering diagram  $\foD=
\{(\fod_i,f_i)\,|\, 1\le i \le n\}$ over 
the ring $R=\kk [[ t_1,\ldots,t_n ]]$ consists of a number of lines
through the origin.
We assume
\begin{enumerate}
\item[$\bullet$]
 $\fod_i=\RR m_i$ with $m_i$ primitive,{\footnote{We do {\em not}
assume the $m_i$ are distinct. Repetition is allowed.}}
\item[$\bullet$]
$t_i$ is the only
power series  variable
occurring in $f_i$, so
$$f_i\in 
\kk[z^{m_i}]\widehat\otimes_{\kk} \kk[[t_i]]\ .$$ 
\end{enumerate}
\end{definition}
\vspace{10pt}

We give here a method of computing 
$\Scatter(\foD)$ to $k$th order for a standard diagram
$\foD$. We work over the ring
\[
R_k=
{\kk [[ t_1,\ldots,t_n]]
\over
(t_1^{k+1},\ldots,t_n^{k+1})}.
\]
Let $\alpha_k:R\rightarrow R_k$ be the projection, and define 
\begin{eqnarray*}
\Scatter_k(\foD)& =& \Scatter(\alpha_k(\foD))\\
& = &\alpha_k(\Scatter(\foD)).
\end{eqnarray*}
We view $\Scatter_k(\foD)$ as $\Scatter(\foD)$ to $k$th order.

A crucial role is played by the ring 
\[
\widetilde R_k={\kk[\{u_{ij}|1\le i\le n,1\le j\le k\}]\over
\langle u_{ij}^2\,|\, 1\le i\le n,1\le j\le k\rangle}.
\]
via the basic inclusion $R_k\hookrightarrow\widetilde R_k$
given by $t_i\mapsto \sum_{j=1}^k u_{ij}$. 
While the $f_i$'s may not factor over $R_k$,
there is a natural factorization of $f_i$ over $\widetilde R_k$
obtained as follows.
We expand $\log f_i$ in a Taylor series in $t_i$, 
\[
\log f_i=\sum_{j=1}^{k}\sum_{w\ge 1}wa_{ijw}z^{wm_i} t^j_i \ \
\in {R}_k,
\]
with $a_{ijw}\in\kk$.
For any given $i$ and $j$, $a_{ijw}=0$ for all but a finite
number of $w$. We then make the substitution 
$t_i=\sum_{j=1}^k u_{ij}$ and expand, getting a sum of monomials
each having a coefficient in $\widetilde R_k$ whose square is zero. 
We find
\[
\log f_i=\sum_{j=1}^k\sum_{\#J=j}\sum_{w\ge 1} j!wa_{ijw}z^{wm_i}
\prod_{l\in J} u_{il}
\]
where $J$ runs over subsets of $\{1,\ldots,k\}$.
Taking the exponential gives us a factorization of $f_i$ as
\[
f_i=\prod_{j=1}^k\prod_{\#J=j}\prod_{w\ge 1} \bigg(1+j!wa_{ijw}z^{wm_i}
\prod_{l\in J} u_{il}\bigg).
\]
We can now apply the  deformation technique to the
 factorization to pull apart the
lines $\fod_i$. We obtain the scattering diagram
\begin{equation*}
\widetilde\foD_k
=\bigg\{\big(\fod_{iJw},1+(\# J)!wa_{i(\# J)w}z^{wm_i}\prod_{l\in J} u_{il}\big)
\,\bigg|\,1\le i\le n, w\ge 1, J\subseteq \{1,\ldots,k\}, \#J\ge 1\bigg\},
\end{equation*}
with $\fod_{iJw}$ a line parallel to $\fod_i$ and chosen generally.

We now describe rather concretely a process for producing 
$\Scatter(\widetilde\foD_k)$. We usually only considered
$\Scatter(\widetilde\foD_k)$ as an equivalence class of scattering diagrams.
However, in Section~2,  we will use the representative
of the equivalence class produced by the procedure given here. 

We produce a sequence of scattering diagrams
\[
\widetilde{\foD}_k=\widetilde\foD^0_k,\ \widetilde\foD^1_k,\ 
\widetilde\foD^2_k,\ \ldots
\]
which eventually stabilizes, and we take $\Scatter(\widetilde\foD_k)
=\widetilde\foD^i_k$
for $i$ sufficiently large. We will assume inductively:
\begin{enumerate}
\item Each ray or line in $\widetilde\foD^i_k$ is of the form
\[
(\fod,1+c_{\fod}u_{I(\fod)}z^{m_{\fod}}),
\]
with $c_{\fod}\in\kk$ and $u_{I(\fod)}=\prod_{(i,j)\in I(\fod)} u_{ij}$,
for some index set 
$$I(\fod)\subseteq \{1,\ldots,n\}\times \{1,\ldots,k\}.$$
\item For each $p\in\Sing(\widetilde\foD^i_k)$, there is no set of rays
$\foD[p]\subseteq\widetilde\foD^i_k$ of cardinality $\ge 3$ such that
$p\in \fod$ for each $\fod\in\foD[p]$ and $I(\fod_1)\cap I(\fod_2)=
\emptyset$ for any two distinct $\fod_1,\fod_2\in\foD[p]$. Note that two
rays $\fod_1,\fod_2$ only produce a non-trivial new ray via 
Lemma~\ref{nilpotentcollision}
if $I(\fod_1)\cap I(\fod_2)=\emptyset$.
\end{enumerate}

Both of these conditions already hold for $\widetilde\foD^0_k$. 
In particular,
for general choice of $\widetilde\foD^0_k$, each singular point
of $\widetilde\foD^0_k$ is contained in only two lines.

To pass from $\widetilde\foD^{i-1}_k$ to $\widetilde\foD^i_k$, we simply look
at every pair of elements $\fod_1,\fod_2\in \widetilde\foD^{i-1}_k$ 
satisfying
\begin{itemize}
\item[(i)] $\{\fod_1,\fod_2\}\not\subseteq \widetilde\foD^{i-2}_k$,
\item[(ii)] $\fod_1\cap\fod_2$ consists of one point $m_0$, which is neither
the boundary of $\fod_1$ nor $\fod_2$,
\item[(iii)] $I(\fod_1)\cap I(\fod_2)=\emptyset$.
\end{itemize}
Writing
\[
f_{\fod_i}=1+c_iu_{I(\fod_i)}z^{w_im_i}
\]
with $m_i$ primitive, we follow Lemma \ref{nilpotentcollision} and set 
\begin{multline}
\label{dd1d2def}
\fod(\fod_1,\fod_2)=\\
\big(m_0+\RR_{\ge 0}(w_1m_1+w_2m_2),1+c_1c_2u_{I(\fod_1)\cup I(\fod_2)}
w_{\out}|m_1\wedge m_2| z^{w_1m_1+w_2m_2}\big).
\end{multline}
We then take
\[
\widetilde\foD^i_k=\widetilde\foD^{i-1}_k\cup\{\fod(\fod_1,\fod_2)|
\hbox{$\fod_1,\fod_2$ satisfies (i)-(iii) above}\}.
\]

Of course $\widetilde\foD^i_k$ satisfies the inductive hypothesis (1), but
 we need to check condition (2). To do so, we introduce some
notation we shall also need later.

For any ray of the form $\fod=\fod(\fod_1,\fod_2)$, define
$\Parents(\fod)=\{\fod_1,\fod_2\}$. If $\fod$ is a line, then
we set $\Parents(\fod)=\emptyset$.  
Define recursively
\[
\Ancestors(\fod)=\{\fod\}\cup\bigcup_{\fod'\in \Parents(\fod)} 
\Ancestors(\fod')
\]
and
\[
\Leaves(\fod)=\{\fod'\in \Ancestors(\fod)| \hbox{$\fod'$ is a line}\}.
\]
If any element $\fod'\in\Leaves(\fod)$ is moved slightly
by replacing $\fod'=m+\RR m_i$ by $m'+\RR m_i$  with
$m-m'$ small and $m-m'\not\in\RR m_i$, then $\fod$ is replaced
by some $m''+\fod$ for $m''$ small and $m''\not\in \RR m_{\fod}$. 
Now 
if there is a point $p\in \Sing(\widetilde\foD^i_k)$ with $\foD[p]\subseteq
\widetilde\foD^i_k$ violating condition~(2), we note  $\Leaves(\fod_1)
\cap \Leaves(\fod_2)=\emptyset$ for all $\fod_1,\fod_2\in \foD[p]$,
since $I(\fod_1)\cap I(\fod_2)=\emptyset$. Thus by deforming
the lines in $\widetilde\foD_k^0$, 
the rays in $\foD[p]$ can be deformed in an independent
manner. Thus if $p$ violates (2), then $\widetilde\foD_k^0$ has not been chosen
generally. Thus for general choice of $\widetilde\foD_k^0$,
the inductive hypothesis (2) holds.

The process is easily seen to terminate as
at each step the sets $I(\fod)$ for new rays $\fod$ increase
in cardinality, and the cardinality of these sets is bounded by $nk$.
Furthermore, since automorphisms associated to rays $\fod_1$, $\fod_2$
with $I(\fod_1)\cap I(\fod_2)\not=\emptyset$ commute, we obtain
 from Lemma \ref{nilpotentcollision} that $\theta_{\gamma_p,
\widetilde\foD^i_k}=\id$ for $i>nk$ and $\gamma_p$ a small loop around
any singular point of $\widetilde\foD^i_k$.

We have proven $\Scatter(\widetilde{\foD}_k) = \widetilde\foD^\infty_k$ 
where stabilization occurs in the superscript for $i> nk$.
Then
$$\Scatter_k(\foD) = (\widetilde\foD_k^\infty)_{\as}.$$
 As we will see in the next section, we will use the
above procedure
to reduce the computation of $\Scatter_k(\foD)$ to a tropical calculation.

\begin{example}
\label{basicexample2}
We apply the above procedure to Example \ref{basicexample} with
$\ell_1=\ell_2=1$ and
 $$f_1=1+t_1x, \ \ \ f_2=1+t_2y\ .$$
For $k=2$, we
write
\[
\log f_1=t_1x-{t_1^2x^2\over 2}+\cdots
=u_{11}x+u_{12}x-u_{11}u_{12}x^2
\]
and factor $f_1$ as
\[
f_1=(1+u_{11}x)(1+u_{12}x)(1-u_{11}u_{12}x^2).
\]
Similarly,
\[
f_2=(1+u_{21}y)(1+u_{22}y)(1-u_{21}u_{22}y^2).
\]
Applying the above procedure gives the complicated 
scattering diagram depicted
in Figure \ref{factordeformexample}. 

In the scattering diagram,
there are four rays of slope $2$, eight of slope $1$, and
four of slope $1/2$. That the asymptotic form
of the diagram agrees with Figure \ref{d1=1,d2=1}
appears miraculous.
To see the agreement, we calculate
 the functions attached to the various
rays in Figure \ref{factordeformexample} using Lemma 
\ref{nilpotentcollision}. For example, the rays
of slope $2$ have attached functions, taken in a clockwise order,
$$1+u_{11}u_{21}u_{22}xy^2, \ 1+u_{12}u_{21}u_{22}xy^2,\ 
1-u_{11}u_{21}u_{22}xy^2, \ 1-u_{12}u_{21}u_{22}xy^2\ .$$
 The product
of these four functions is $1$, so the ray in the asymptotic diagram
associated to this diagram has the attached function $1$, and hence can
be omitted. 
The same happens to the rays of slope $1/2$. For the rays
of slope $1$, again in clockwise order, again using 
Lemma \ref{nilpotentcollision}, we have the attached functions
$$1+u_{11}u_{21}xy,\ 1+u_{11}u_{22}xy, \ 1+4u_{11}u_{12}u_{21}u_{22}x^2y^2,\ 
1+u_{12}u_{21}xy,\ 
1+u_{12}u_{22}xy,$$ 
$$1-4u_{11}u_{12}u_{21}u_{22}x^2y^2,\ 
1-4u_{11}u_{12}u_{21}u_{22}x^2y^2,\ 1+2u_{11}u_{12}u_{21}u_{22}x^2y^2\ .$$
Taking the product of these eight functions gives
$1+(u_{11}+u_{12})(u_{21}+u_{22})xy$. All the terms involving $x^2y^2$ cancel.

\begin{figure}
\input{factordeformexample.pstex_t}
\caption{}
\label{factordeformexample}
\end{figure}

\end{example}


\section{Tropical curves} \label{tcc}
\subsection{Definitions}
We relate here the calculations involving scattering diagrams
in Section \ref{asymfac} to tropical curve counts. We first recall the
definition of a parameterized tropical curve in $M_{\RR}$ from \cite{Mk}.
 Let $\overline{\Gamma}$ be a weighted, connected
finite graph without divalent vertices. Denote the set of vertices and
edges by $\overline{\Gamma}^{[0]}$ and $\overline{\Gamma}^{[1]}$
respectively, and let
$$w_{\overline{\Gamma}}:\overline{\Gamma}^{[1]}
\rightarrow \NN\setminus \{0\}$$ be the weight function. 
An edge
$E\in \overline{\Gamma}^{[1]}$ has adjacent vertices $\partial E=\{V_1,V_2\}$.
Let $\overline{\Gamma}_{\infty}^{[0]}\subseteq \overline{\Gamma}^{[0]}$
be the set of 1-valent vertices. We set
\[
\Gamma=\overline{\Gamma}\setminus\overline{\Gamma}_{\infty}^{[0]}
\]
Denote the set of vertices and edges of $\Gamma$ as $\Gamma^{[0]}$,
$\Gamma^{[1]}$, and let $$w_{\Gamma}:\Gamma^{[1]}
\rightarrow \NN\setminus\{0\}$$
be the weight function.
Some edges of $\Gamma$ are now non-compact --- these are called
\emph{unbounded edges}. Let $\Gamma_{\infty}^{[1]}\subseteq
\Gamma^{[1]}$ be the set of unbounded edges.

\begin{definition}
\label{tropcurve}
A \emph{parameterized tropical curve} in $M_{\RR}$ is a proper map
$h:\Gamma\rightarrow M_{\RR}$ satisfying the following conditions.
\begin{enumerate}
\item For every edge $E\subseteq\Gamma$, the restriction $h|_E$ is
an embedding with image $h(E)$ contained in an affine line with rational
slope.
\item For every vertex $V\in\Gamma$, the following \emph{balancing
condition} holds. Let $$E_1,\ldots,E_m\in \Gamma^{[1]}$$
 be the edges
adjacent to $V$, and let $m_i\in M$ be the primitive integral vector
emanating from $h(V)$ in the direction of $h(E_i)$. Then
\[
\sum_{j=1}^m w_\Gamma(E_j)m_j=0.
\]
\end{enumerate}
An \emph{isomorphism} of parameterized
tropical curves $h:\Gamma\rightarrow M_{\RR}$
and $h':\Gamma'\rightarrow M_{\RR}$ is a homeomorphism $\Phi:\Gamma
\rightarrow\Gamma'$ respecting the weights of the edges and such that
$h=h'\circ\Phi$. A \emph{tropical curve} is an isomorphism class
of parameterized tropical curves. The \emph{genus} of a tropical curve
$h:\Gamma\rightarrow M_{\RR}$ is the first Betti number of $\Gamma$.
A \emph{rational} tropical curve is a tropical curve of genus zero.
\end{definition}

\begin{definition}
\label{multiplicitydefinition}
Let $h:\Gamma\rightarrow M_{\RR}$ be a tropical
curve such that $\bar\Gamma$ only has vertices of valency one and
three. The \emph{multiplicity} of a vertex $V\in\Gamma^{[0]}$ in $h$ is
\[
\Mult_V(h)=w_1w_2|m_1\wedge m_2|=w_1w_3|m_1\wedge m_3|=w_2w_3|m_2\wedge 
m_3|,
\]
where $E_1,E_2,E_3\in\Gamma^{[1]}$ are the edges containing $V$,
$w_i=w_{\Gamma}(E_i)$, and $m_i\in M$ is a primitive integral vector
emanating from $h(V)$ in the direction of $h(E_i)$. The equality of the
three expressions follows from the balancing condition.
\end{definition}

\begin{definition}
\label{multiplicitydefinition2}
The \emph{multiplicity} of the tropical curve $h$ is 
\[
\Mult(h)=\prod_{V\in\Gamma^{[0]}}\Mult_V(h).
\]
\end{definition}

\subsection{Scattering diagrams}

Let $\foD=\{(\RR m_i,f_i)\}$ be a standard scattering
diagram of Definition \ref{ssss}.
We follow the notation of Section \ref{asymfac}
with scattering diagrams $\widetilde
\foD_k$ and 
$$\Scatter(\widetilde\foD_k)= \widetilde\foD_k^\infty$$ 
constructed  from $\foD$.

\begin{theorem}
\label{onetoonecorrespondence}
For general choice of $\widetilde\foD_k$,
there is a bijective correspondence between  elements
$(\fod,f_{\fod})\in \Scatter(\widetilde\foD_k)$ and rational tropical curves
$h:\Gamma\rightarrow M_{\RR}$ with the
following properties:
\begin{enumerate}
\item[(i)] There is an edge $E_{\out}\in\Gamma^{[1]}_{\infty}$ with
$h(E_{\out})=\fod$.
\item[(ii)] If $E\in \Gamma^{[1]}_{\infty}\setminus
\{E_{\out}\}$ or if $E_{\out}$ is the only edge of $\Gamma$
and $E=E_{\out}$, 
then $h(E)$ is contained in some $\fod_{iJw}$ where 
$$1\le i\le n,\ \ \ J\subseteq\{1,\ldots,k\}, \ \ \ w\ge 1.$$ Furthermore,
if $E\not=E_{\out}$, 
the unbounded direction of $h(E)$ is given by $-m_i$.
\item[(iii)]
 If $E,E'\in\Gamma^{[1]}_{\infty}\setminus\{E_{\out}\}$
and $h(E)\subseteq \fod_{iJw}$ and $h(E')\subseteq
\fod_{iJ'w'}$, then $J\cap J'=\emptyset$.
\item[(iv)] If $E\in\Gamma^{[1]}_{\infty}\setminus\{E_{\out}\}$
or if $E_{\out}$ is the only edge of $\Gamma$ and $E=E_{\out}$,
and $h(E)\subseteq \fod_{iJw}$, we have $w_{\Gamma}(E)=w$.
\end{enumerate}
Furthermore, if $\fod$ is a ray, the corresponding curve $h$
is trivalent and
\begin{equation}
\label{ffodform}
f_{\fod}=1+w_{\out}\Mult(h)\prod_{i,J,w}\bigg((\#J)!
a_{i(\#J)w}\prod_{j\in J}u_{ij}\bigg)z^{m_{\out}},
\end{equation}
where the $i,J,w$ run over all indices for which  $\fod_{iJw}\in\Leaves(\fod)$,
\[
m_{\out}=\sum_{i,J,w} wm_i,
\]
and $m_{\out}=w_{\out}m'_{\out}$ for $m'_{\out}\in M$ primitive and
$w_{\out}$ the index of $m_{\out}$.
\end{theorem}

\proof The bijective correspondence is constructed as follows.
First, the elements of $\widetilde\foD_k$ are all lines and 
are clearly in one-to-one
correspondence with tropical curves $h:\Gamma\rightarrow M_{\RR}$
in which $\Gamma$ consists of just one edge satisfying properties
(i)-(iv). So we just need to worry about elements of $\Scatter(\widetilde
\foD_k)\setminus\widetilde\foD_k$, all of which are rays.

Given a ray $\fod\in\widetilde\foD_k^\infty$,
define the graph $\Gamma_{\fod}$ by 
\begin{eqnarray*}
\Gamma_{\fod}^{[0]}&=&\{\fod'\,|\,\hbox{$\fod'\in\Ancestors(\fod)$
and $\fod'$ a ray}\}\\
\Gamma_{\fod}^{[1]}&=&\Ancestors(\fod)
\end{eqnarray*}
as abstract sets, writing the edge corresponding to $\fod'$ as
$E_{\fod'}$ and the vertex $V_{\fod'}$. 
If $\fod'\in\Ancestors(\fod)\setminus\{\fod\}$, then $\fod'$ is parent
to a unique ray in $\Ancestors(\fod)$, which we write as $\Child(\fod')$.
If $\fod$ is itself a line, then we are in a degenerate case, and
$\Gamma_{\fod}$ 
is just a line (an edge with no vertices). If $\fod$ is a ray,
then to define $\Gamma_{\fod}$, we specify that for $\fod'\in\Ancestors(\fod)$,
either 
\begin{itemize}
\item $\fod'\not=\fod$ and $\fod'$ is a ray, in which case the vertices of
$E_{\fod'}$ are $V_{\fod'}$ and $V_{\Child(\fod')}$.
\item $\fod'=\fod$, in which case $E_{\fod'}$ is an unbounded edge
with vertex $V_{\fod}$.
\item $\fod'$ is a line, in which case $E_{\fod'}$ is an unbounded edge
with vertex $V_{\Child(\fod')}$.
\end{itemize}
We next define the weight function $w_{\Gamma_{\fod}}$ as follows.
For any $\fod'\in\Ancestors(\fod)$, $f_{\fod'}$ takes the form
$1+c_{\fod'}z^{m_{\fod'}}$ for some $c_{\fod'}\in \tilde R_k$ and 
$m_{\fod'}\in M\setminus\{0\}$. We
define $w_{\Gamma_{\fod}}(E_{\fod'})$ to be the index of $m_{\fod'}$, with 
$m'_{\fod'}$ given by
\[
m_{\fod'}=w_{\Gamma_{\fod}}(E_{\fod'})m'_{\fod'}.
\]

Finally, we define $h$ in the obvious way, mapping $E_{\fod'}$ in
the respective cases to 
the line segment joining $\Init(\fod')$ and $\Init(\Child(\fod'))$, 
 $\fod$, and  the ray $\RR_{\le 0} m_{\fod'}+\Init(\Child(\fod'))$.
 
Since $\Gamma$ is a tree, the genus 0 condition is clearly satisfied.
We need to check the balancing condition at every vertex $V=V_{\fod'}$.
Let $\Parents(\fod')=\{\fod_1,\fod_2\}$. From 
\eqref{dd1d2def}, we see  
\[
m_{\fod'}=m_{\fod_1}+m_{\fod_2},
\]
or equivalently
\[
w_{\Gamma_{\fod}}(E_{\fod'})m'_{\fod'}=
w_{\Gamma_{\fod}}(E_{\fod_1})m'_{\fod_1}+
w_{\Gamma_{\fod}}(E_{\fod_2})m'_{\fod_2},
\]
which is the balancing condition at $V$, keeping in mind that
$m'_{\fod_1}$ and $m'_{\fod_2}$ point \emph{towards} $V$ and
$m'_{\fod'}$ points \emph{away} from $V$.
By the generality condition on $\widetilde\foD_k$,
 $\Gamma_{\fod}$ has at most trivalent vertices.

Next, we check the expression given for $f_{\fod}$ inductively.
If $\fod$  is a line, 
then $\Gamma$ is just a line,
and $\Mult(h)=1$ since there are no trivalent vertices whatsoever.
Formula \eqref{ffodform} is correct by the definition of the
 original deformed
scattering diagram $\widetilde\foD_k$ in Section \ref{asymfac}. 

Suppose $\fod$ is
a ray and \eqref{ffodform} holds for both
$\Parents(\fod)=\{\fod_1,\fod_2\}$. Let $h_1$ and $h_2$ be the 
tropical curves corresponding to the respective
parents. By \eqref{dd1d2def},
\begin{multline*}
f_{\fod}=
1+ 
w_{\Gamma_{\fod}}(E_{\fod_1})
\Mult(h_1)w_{\Gamma_{\fod}}(E_{\fod_2})
\Mult(h_2)\\ \cdot w_{\out}|m'_{\fod_1}\wedge m'_{\fod_2}|\prod_{i,J,w}\bigg( (\#J)!a_{i(\#J)w}
\prod_{j\in J} u_{ij}\bigg)z^{m_{\out}}
\end{multline*}
where the first product is over all indices $i,J,w$ for which
$\fod_{iJw}\in\Leaves(\fod)$.
The coefficient in front of the product is just
\[
\Mult(h_1)\Mult(h_2)\Mult_{V_{\fod}}(h)w_{\out}=\Mult(h)w_{\out}.
\]
The derivation of \eqref{ffodform} is complete.

Finally, by unravelling the definitions,
each rational tropical curve satisfying conditions (i)-(iv)
is easily seen 
to correspond to a unique element of $\widetilde\foD_k^\infty$.
\qed

\bigskip

As we shall see, Theorem \ref{onetoonecorrespondence} implies that
the scattering
diagrams are computing a tropical enumerative invariant.
We shall relate the tropical enumerations which arise
first to holomorphic enumerations
and then to relative Gromov-Witten invariants.

\begin{example}
\label{basicexample3}
Returning to the situation of Example \ref{basicexample2},
consider the ray of slope $1$ which is the third one from the upper
left in Figure \ref{factordeformexample}. Figure \ref{curvefigure}
shows the corresponding tropical curve with the outgoing edge
of weight two labelled. The multiplicity of the curve is
$2$, and the function attached to the outgoing ray is $1+4u_{11}u_{12}
u_{21}u_{22}x^2y^2$.
\begin{figure}
\input{curvefigure.pstex_t}
\caption{}
\label{curvefigure}
\end{figure}
\end{example}

\subsection{Tropical counts}
\label{tcounts}
Let $m_1,\ldots,m_n\in M$ be primitive elements, and let
$${\bf m}=(m_1, \ldots,m_n)$$ denote the $n$-tuple.
Let the  lines 
$$\fod_{ij}=m_{ij}+\RR m_i, \ \ \  
 m_{ij}\in M_{\RR}$$ be chosen
generally for
$1\le i\le n$ and $1\le j\le l_i$.

Let  
${\bf w}_i
=(w_{i1},\ldots,w_{il_i})$ be 
\emph{weight vectors}  with
$$0< w_{i1}\le w_{i2}\le\cdots \le w_{il_i}, \ \ \
w_{ij} \in \mathbb{Z}\ .$$
The weight vector ${\bf w}_i$ has
{\em length} $l_i$ and
 {\em size} 
\[
|{\bf w}_i|=\sum_{j=1}^{l_i} w_{ij}.
\]
Denote the $n$-tuple of weight vectors by
${\bf w}= ( {\bf w}_1, \ldots, {\bf w}_n).$
We only consider weight vectors for which
$$m_{\out} = \sum_{i=1}^n |{\bf w}_i|  m_{i}\neq 0 .$$

\begin{definition}\label{def.Ntrop}
Let 
$
N^{\trop}_{\bf m}(\bf{w})
$
be
the number of rational tropical curves $h:\Gamma\rightarrow M_{\RR}$
counted with the 
multiplicity of Definition \ref{multiplicitydefinition}
where
\[
\Gamma^{[1]}_{\infty}=\{E_{ij}|1\le i\le n, 1\le j\le l_i\} \cup \{E_{\out}\}
\]
with
\begin{itemize}
\item
$h(E_{ij})\subseteq\fod_{ij}$ and $-m_i$ pointing in the unbounded
direction of $h(E_{ij})$, 
\item  $w_{\Gamma}(E_{ij})=w_{ij}$, 
\item   $h(E_\out)$ pointing in the direction of
$m_{\out}$.
\end{itemize}
Note that $w_{\Gamma}(E_{\out})$, the weight of $E_{\out}$, coincides
with the index of $m_{\out}$
\end{definition}

\bigskip

\begin{proposition}
The numbers $N^{\trop}_{\bf m}({\bf w})$ do not depend on the (generic)
choice
of the vectors $m_{ij}$.
\end{proposition}

\proof The result follows
from  standard tropical arguments --- for example the
method of \cite{GaMa} suffices. Alternatively, we will later show
that these tropical invariants
agree with holomorphic counts which are independent of any choices.
\qed

\bigskip

Let $\foD=
\{(\fod_i,f_i)\,|\, 1\le i \le n\}$ be a standard scattering
diagram over 
the ring $\kk [[ t_1,\ldots,t_n ]]$ consisting of a number of lines
through the origin. Let $\fod_i=\RR m_i$ with $m_i$ primitive.
Since $t_i$ is the only
power series  variable
occurring in $f_i$,
we can write the logarithm as
\[
\log f_i=\sum_{j=1}^{\infty}\sum_{w\ge 1}wa_{ijw}z^{wm_i} t^j_i \ \
\]
with $a_{ijw}\in\kk$.

Let $\Scatter(\foD)$ be the associated scattering diagram.
We 
assume $\Scatter(\foD)$ has at most one ray in any given direction.

\begin{theorem} \label{mtr}
Let $(\fod,f_{\fod})\in \Scatter(\foD)\setminus\foD$ be a ray. Then,
\begin{equation}
\label{maintropical}
\log f_{\fod}=\sum_{\bf w} \sum_{\bf k} \frac{
w_{\out}({\bf w})N^{\trop}_{\bf m}({\bf w})}{
 |\Aut({\bf w},{\bf k})|} 
\bigg(\prod_{1\le i\le n\atop 1\le j\le l_i}
a_{ik_{ij}w_{ij}} t_i^{k_{ij}}\bigg)
z^{\sum_i |{\bf w}_i|m_i}
\end{equation}
where 
\begin{itemize}
\item
The first sum is over all $n$-tuples of weight vectors ${\bf w}=
({\bf w}_1, \ldots, {\bf w}_n)$
satisfying
 $0\neq \sum_i |{\bf w}_i|m_i\in\fod$. 
Let $l_i=\length({\bf w}_i)$.
\item
The second sum is over all $n$-tuples of vectors 
${\bf k}=
({\bf k}_1, \ldots, {\bf k}_n)$ where
$${\bf k}_i=(k_{i1},\ldots,k_{il_i})$$
for positive integers 
$k_{ij}$ satisfying
$$k_{ij}\le k_{i(j+1)}\ \text{ if }\
w_{ij}=w_{i(j+1)}.$$
Note $\length({\bf k}_i)=\length({\bf w}_i)$.
\item $\sum_i |{\bf w}_i|m_i=
w_{\out}({\bf w})m'_{\out}$ with $m'_{\out}\in M$ primitive.
\item 
$\Aut({\bf w}_i,{\bf k}_i)$ denotes the subgroup of the permutation
group
$\Sigma_{l_i}$
stabilizing the data
$$\left( (w_{i1},k_{i1}), \ldots, (w_{il_i},k_{il_i}) \right),$$
$|\Aut({\bf w}_i,{\bf k}_i)|$
denotes the order,
and $$|\Aut({\bf w}, {\bf k})| =\prod_{i=1}^n |\Aut({\bf w}_i,{\bf k}_i)|\ .$$ 
\end{itemize}
\end{theorem}

\proof We will prove \eqref{maintropical} modulo the ideal
$I_k=(t_1^{k+1},\ldots,t_n^{k+1})$ and let $k\rightarrow\infty$.
Modulo $I_k$, the diagram
$\Scatter(\foD)$ is equivalent to $(\widetilde\foD_k^\infty)_{\as}$.
Thus modulo $I_k$,
\begin{equation}
\label{logsum}
\log f_{\fod}=\sum_{\fod'} \log f_{\fod'},
\end{equation}
where the sum is over rays $\fod'\in\widetilde\foD^\infty_k$ 
parallel to $\fod$.

Theorem \ref{onetoonecorrespondence} gives a correspondence
between such rays and tropical
curves $$h:\Gamma\rightarrow M_{\RR}.$$
 Consider one such ray $\fod'$ and the
corresponding tropical curve $h:\Gamma\rightarrow M_{\RR}$. Then,
there are weight vectors
${\bf w}_1(h),\ldots,{\bf w}_n(h)$ of lengths $l_1,\ldots,l_n$, and sets
$$J_{ij}(h)\subseteq \{1,\ldots,k\},\ \ \ 1\le j\le l_i$$
 with $J_{i1}(h),
\ldots,J_{il_i}(h)$ pairwise disjoint, such that
we can write 
\[
\Leaves(\fod')=\{\fod_{iJ_{ij}(h)w_{ij}(h)}\,|\,1\le i\le n, 1\le j\le l_i\}
\]
and
\[
\Gamma^{[1]}_{\infty}=\{E_{ij}|1\le i\le n, 1\le j\le l_i\} \cup \{E_{\out}\}
\]
with 
$h(E_{ij})\subseteq \fod_{iJ_{ij}(h)w_{ij}(h)}$.
Moreover, the contribution from $\fod'$ to \eqref{logsum} is
\[
\log f_{\fod'}=w_{\Gamma}(E_{\out})\Mult(h)
\prod_{1\le i\le n
\atop 1\le j\le l_i} \big((\# J_{ij}(h))!\ a_{i(\#J_{ij}(h))
w_{ij}(h)}\prod_{l\in J_{ij}(h)}
u_{il}\big) z^{\sum_i |{\bf w}_i(h)|m_i}.
\]
The contribution from all tropical curves $h$ giving rise to
the same weight vectors ${\bf w}_i$ and the same sets $J_{ij}$ with
$k_{ij}=\#J_{ij}$ defining vectors ${\bf k}_1,\ldots,{\bf k}_n$, is then
\[
w_{\out}({\bf w})N^{\trop}_{\bf m}({\bf w})
\prod_{1\le i\le n
\atop 1\le j\le l_i} \big(k_{ij}!\ a_{ik_{ij}w_{ij}}\prod_{l\in J_{ij}}
u_{il}\big) z^{\sum_i |{\bf w}_i|m_i}.
\]

We first 
 keep both ${\bf w}$ and ${\bf k}$ fixed  and choose sets
$J_1,\ldots,J_n\subseteq \{1,\ldots,k\}$ with $\#J_i=|{\bf k}_i|$. We then
sum over all possible ways of writing $J_i$ as a disjoint union
of sets $J_{ij}$ with $\#J_{ij}=k_{ij}$. There are $|{\bf k}_i|!/\prod_j
k_{ij}!$ ways of writing $J_i$ as such a disjoint union. However,
we have
overcounted curves. If $\sigma\in\Aut({\bf w}_i,{\bf k}_i)$,
we have
\[
\{(i,J_{ij},w_{ij})| 1\le j\le \length({\bf w}_i)\}=
\{(i,J_{i\sigma(j)},w_{i\sigma(j)})| 1\le j\le \length({\bf w}_i)\}.
\]
Thus the contribution from a choice of the sets $J_i$ is 
\[
{w_{\out}({\bf w})N^{\trop}_{\bf m}({\bf w})
\over |\Aut({\bf w},{\bf k})|} 
\bigg(\prod_{1\le i\le n}
\big(\prod_{1\le j\le l_i}
k_{ij}!\ a_{ik_{ij}w_{ij}}\big)
{|{\bf k}_i|!\over \prod_j k_{ij}!}
\prod_{l\in J_i}
u_{il} \bigg)
z^{\sum_i |{\bf w}_i|m_i}.
\]
Summing over possible $J_i$'s, we get
\[
{w_{\out}({\bf w})N^{\trop}_{\bf m}({\bf w})
\over |\Aut({\bf w},{\bf k})|} 
\bigg(\prod_{1\le i\le n}
\big(\prod_{1\le j\le l_i}
a_{ik_{ij}w_{ij}}\big)\sum_{J_i\atop \#J_i=|{\bf k}_i|}
|{\bf k}_i|! \prod_{l\in J_i}
u_{il} \bigg)
z^{\sum_i |{\bf w}_i|m_i}.
\]
But 
\begin{eqnarray*}
t_i^{|{\bf k}_i|}&=&(u_1+\cdots+u_k)^{|{\bf k}_i|}\\
&=&|{\bf k}_i|!\sum_{J_i\atop \#J_i=|{\bf k}_i|}\prod_{l\in J_i} u_{il},
\end{eqnarray*} from which we find the total contribution from curves
with given ${\bf w}_i$ and ${\bf k}_i$ is
\[
{w_{\out}({\bf w})N^{\trop}_{\bf m}({\bf w})
\over |\Aut({\bf w},{\bf k})|} 
\bigg(\prod_{1\le i\le n
\atop 1\le j\le l_i}
a_{ik_{ij}w_{ij}}t_i^{k_{ij}}\bigg)
z^{\sum_i |{\bf w}_i|m_i},
\]
giving the desired result.
\qed


\section{From tropical to holomorphic counts}
\label{section-tropicaltoholomorphic}

\subsection{Holomorphic counts}
\label{hcounts}
We now describe the holomorphic analogue
of $N^{\trop}_{\bf m}({\bf w})$. 
Following the notation of Section \ref{tcounts},
let $${\bf m}=(m_1,\ldots,m_n)$$ be an $n$-tuple of
  primitive vectors of $M$.
Let ${\bf w}=({\bf w}_1,\ldots, {\bf w}_n)$ be an $n$-tuple of weight vectors
$${\bf w}_i=(w_{i1},\ldots,w_{il_i})$$
 with
$$0\neq m_{\out}
=\sum_{i=1}^n |{\bf w}_i| m_i, \ \  m_{\out}=w_{\out}m'_{\out}$$
for $m'_{\out}\in M$ primitive.
 
To match the conventions 
set for standard scattering diagrams in Definition \ref{ssss},
we do not require the $m_1, \ldots, m_n$ to
be distinct. However, we will only treat the 
distinct case since the multiplicities here will not matter
(and complicate the notation).
A more subtle issue concerns $m_\out$.
Either the ray generated by $m_\out$ is distinct
from the  rays  generated by  
$-m_1,\ldots,-m_n$ or not. We will present
a treatment of  the former case. The discussion in 
the degenerate case where the ray generated by $m_\out$
coincides with a ray generated by $-m_k$ 
 is almost identical. We state the results there
and leave the details to the reader.

So we assume the ray generated by $m_\out$ is distinct from
the rays generated by the $-m_i$.
Let $\Sigma$
denote the complete rational fan in $M_{\RR}$ whose
rays are 
generated
by 
\begin{equation*}
-m_1,\ldots,-m_n,m_{\out}. 
\end{equation*}
Let $X$
denote the corresponding toric surface over $\CC$.
Let $$D_1,\ldots,D_n,D_{\out} \subseteq X$$
be the toric divisors corresponding to the given rays. 
Let $X^o$ be the complement of the 0-dimensional torus orbits
in $X$, and let
$$D_i^o=D_i\cap X^o, \ \ D_{\out}^o=D_{\out}\cap X^o.$$

Let $\beta_{\bf w}
\in H_2(X,\ZZ)$ be the homology class
defined by the conditions
$$D_i\cdot\beta_{\bf w}=|{\bf w}_i|,
\ \ \ D_{\out}\cdot \beta_{\bf w}
=w_{\out}.$$
Define the open subspace
\[
\foM(X^o,{\bf w})\subseteq 
\overline{\foM}_{0,1+\sum_{i=1}^n l_i}
(X,\beta_{\bf w})
\]
of the moduli space of genus 0 stable maps
represented by maps of the form 
\[
\varphi:(\PP^1,Q_{11},\ldots,Q_{1l_1},\ldots,
Q_{n1},\ldots,Q_{nl_n}, Q)
\rightarrow X
\]
satisfying the following properties:
\begin{enumerate}
\item[(i)]  
 $\varphi(Q_{ij})\in D_i^o$  
and  $\varphi^*(t_i)$ has a zero
of order $w_{ij}$ at $Q_{ij}$ for a local defining
equation $t_{i}$ for $D_i$.
\item[(ii)] 
$\varphi(Q)\in D^o_{\out}$ and $\varphi^*(t_{\out})$ has a zero
of order $w_{\out}$ at $Q$ for a local defining
equation $t_{\out}$ for $D_\out$.
\end{enumerate}

Since all the intersections of $\varphi$ with the
toric divisors of $X$ are accounted for, there is a
factorization
$$\varphi: \mathbb{P}^1 \rightarrow X^o \subseteq X.$$
We may equivalently view $\foM(X^o,{\bf w})$  as an open
subspace of the moduli space of genus 0 stable maps
to $X^o$ of class $\beta_{\bf w}$
relative to the divisors $D_1^o, \ldots,D_n^o, D_\out$
with relative conditions specified by ${\bf w}$ and $w_\out$.
The relative perspective will be pursued in Sections \ref{gwt} and \ref{ftvg}.

\begin{definition}
\label{holomorphicinvariant}

There is an evaluation map
\begin{equation} \label{scvv}
\ev:\foM(X^o,{\bf w})\rightarrow \prod_{i=1}^n (D_i^o)^{l_i}
\end{equation}
If the evaluation map is generically finite, then we define
$N^{\text{hol}}_{\bf m}({\bf w})$ to be the degree. 
Otherwise we define $N^{\text{hol}}_{\bf m}({\bf w})$ to be
0.
\end{definition}

The count $N_{\bf m}^{\text{hol}}({\bf w})$ is invariant 
under refinement of the fan $\Sigma$. The curves
which appear in $\foM(X^o,{\bf w})$ are disjoint from the 
0-dimensional torus orbits of $X$, and hence 
lift to any toric blow-up of $X$. Conversely, any
curve on a toric blow-up of $X$ disjoint from the
exceptional divisors is the lifting of a curve on $X$.
As a consequence, we can always assume that $X$ is
nonsingular by adding additional vectors $m_i$ with
associated weight vectors of length $0$.

By the following result, the dimensions of the domain and
the target of the evaluation map \eqref{scvv}
coincide.

\begin{proposition}
\label{Mdimension}
$\foM(X^o,{\bf w})$ is a nonsingular Deligne-Mumford stack of dimension
$\sum_{i=1}^n \length({\bf w}_i)$.
\end{proposition}

\proof 
We first review the
standard quotient presentation of $X$, which we assume to be
non-singular.
Let $T$ be the free abelian group generated by $\Sigma^{[1]}$,
the set of rays of $\Sigma$. Let $$T_{\CC}=
T\otimes_{\ZZ}\CC.$$
 We may view $T_{\CC}$ as
an affine space,  
\[
T_{\CC}=\Spec \CC[\{z^{\rho}\,|\,\rho\in\Sigma^{[1]}\}].
\]
Given a subset $S \subseteq \Sigma^{[1]}$, define ${\bf A}(S)$
to be the subspace of $T_{\CC}$ determined by the equations
$\{z^{\rho}=0\,|\,\rho\in S\}$. Set
\[
Z=\bigcup_S {\bf A}(S)
\]
where the union is over all subsets $S\subseteq \Sigma^{[1]}$
which are \emph{not} the rays of some cone in $\Sigma$.
Set
\[
U=T_{\CC}\setminus Z.
\]
There is a natural map $r:T\rightarrow M$ taking a generator
$\rho\in \Sigma^{[1]}$ of $T$ to $m_{\rho}\in M$, the primitive
generator of the ray $M$. This yields an exact sequence 
\[
0\mapright{} K\mapright{} T\mapright{r} M\mapright{} 0.
\]
The algebraic torus $K\otimes_{\ZZ} \CC^*$ acts
on $T_{\CC}$, and we have
\[
X=U/ K\otimes_{\ZZ}\CC^*.
\]
Because $X$ is nonsingular, the $K\otimes_{\ZZ} \CC^*$ action
on $U$ is free.

We write the coordinates on $T_{\CC}$
as $z_1,\ldots,z_n,z_{\out}$ corresponding to the
rays generated by $-m_1,\ldots,-m_n$ and $m_{\out}$. These are
our homogeneous coordinates on $X$. The divisor
$D_i$ is given by $z_i=0$ and the divisor $D_{\out}$ is given by
$z_{\out}=0$.

A parameterized map $\varphi:\PP^1\rightarrow
X$ with image of class $\beta$ yields (up to an action of $K\otimes \CC^*)$  
homogeneous
polynomials $\varphi_1,\ldots,\varphi_n,\varphi_{\out}$ in the
coordinate ring $\CC[u,v]$ of $\PP^1$
with  $\deg\varphi_i=|{\bf w}_i|$ which avoid $Z$.
Conversely, such polynomials $\varphi=(\varphi_1,\ldots,
\varphi_n,\varphi_{\out})$
determine a unique parameterized map.

If $\varphi\in \foM(X^o,{\bf w})$, we can always apply an element of 
$\PGL(\PP^1)$
so that $$Q=(1:0)\in\PP^1.$$
 Then, we must have
\begin{enumerate}
\item[(i)] $\varphi_i=c_i \prod_{j=1}^{l_i} (u-a_{ij}v)^{w_{ij}}$ for
some $c_i\in \CC^{*},a_{ij}\in\CC$.
\item[(ii)] $\varphi_{\out}=c_{\out}v^{w_{\out}}$, $c_{\out}\in\CC^{*}$.
\end{enumerate}
Furthermore, such a choice of $\varphi_i,\varphi_{\out}$ gives rise to
$\varphi\in\foM(X^o,{\bf w})$ if and only if all the $a_{ij}$'s
are distinct. The space of choices of $\varphi_i,\varphi_{\out}$
is therefore an open subset of $$(\CC^*)^{n+1}\times
\prod_{i=1}^n\CC^{\length({\bf w}_i)}.$$
 We then have to divide out by
the action of $K\otimes\CC^{*}$, which is of dimension $(n+1)-2$,
as well as the 2-dimensional subgroup of $\PGL(\PP^1)$ keeping
$(1:0)$ fixed (to remove the parameterization). Since these actions have at most finite
stabilizers, we obtain a nonsingular Deligne-Mumford stack of dimension
$\sum_{i=1}^n\length({\bf w}_i)$. 
\qed
\bigskip

We will need later a dimension bound for the same
toric geometry $X^o$, but with different curve classes.
Let ${\beta}_{\widehat{\bf w}}\in H_2(X^o,\mathbb{Z})$ be a 
class defined by the conditions
$$D_i\cdot {\beta}_{\widehat{\bf w}}
 =|\widehat{\mathbf{w}}_i|,
\ \ \ D_{\out}\cdot \beta_{\widehat{\bf w}}
= 0$$
for  weight vectors $\widehat{\bf w}= (\widehat{\mathbf{w}}_1, \ldots
\widehat{\mathbf{w}}_n)$ satisfying
$$\sum_{i=1}^n |\widehat{\mathbf{w}}_i| m_i =0.$$
Then, the proof of Proposition \ref{Mdimension} immediately yields
the following result.

\begin{proposition}
\label{Mdimension2}
$\foM(X^o,\widehat{\bf w})$ is a nonsingular Deligne-Mumford stack
 of dimension
$-1+\sum_{i=1}^n \length(\widehat{{\bf w}}_i)$.
\end{proposition}

\subsection{Equivalence}

The following result relating the holomorphic and tropical
counts will play a crucial role for us.

\begin{theorem}\label{N=Ntrop}
$N_{\bf m}^{{\trop}}({\bf w})=N^{\text{\em hol}}_{\bf m}({\bf w}) \cdot
\prod_{i=1}^n \prod_{j=1}^{l_i} w_{ij}$ .
\end{theorem}

Tropical counting of holomorphic curves on
toric varieties with incidence conditions occurs  in \cite{Mk,nisi}. 
The latter reference already includes higher order
tangency conditions with the toric boundary. The only difference between
\cite{nisi} and
the present situation is the treatment of point conditions ---
here we impose point conditions 
on the toric boundary instead of in the big cell.
Very little has be to added to \cite{nisi} to obtain Theorem \ref{N=Ntrop}.
 The required 
modifications are straightforward. A more detailed discussion is
presented in the Appendix.

\subsection{Degenerate case} \label{degcas}
Suppose the ray generated by $m_\out$
coincides with the ray generated by $-m_k$.
Let $\Sigma$
denote the complete rational fan in $M_{\RR}$ whose
rays are 
generated
by 
\begin{equation*}
-m_1,\ldots,-m_n\ .
\end{equation*}
Let $X$
denote the corresponding toric surface over $\CC$.
Let $$D_1,\ldots,D_n \subseteq X$$
be the toric divisors corresponding to the given rays
and let $D_\out = D_k$.

Let $\beta_{\bf w}
\in H_2(X,\ZZ)$ be the homology class
defined by the conditions
$$D_i\cdot\beta_{\bf w}=|{\bf w}_i| \ \text{ for } i \neq k, \ \
\ \ \ D_{k}\cdot \beta_{\bf w}
=|{\bf w}_k|+w_{\out}.$$
We define moduli spaces $\foM(X^o,{\bf w})$ of maps
of class $\beta_{\bf w}$,
evaluation maps, and invariants $N^{\text{hol}}_{\bf m} ({\bf w})$
in exactly the same manner as before. The only difference
is the contact point $Q$ of
\[
\varphi:(\PP^1,Q_{11},\ldots,Q_{1l_1},\ldots,
Q_{n1},\ldots,Q_{nl_n}, Q)
\rightarrow X
\]
shares a divisor
with $Q_{k1}, \ldots, Q_{kl_k}$.
The rest is straightforward. All the parallel results
hold including Theorem \ref{N=Ntrop}.


\section{Gromov-Witten theory}
\label{gwt}

\subsection{Overview}
We will now  connect the
holomorphic counts $N^{\text{hol}}_{\bf m}({\bf w})$ defined in 
Definition \ref{holomorphicinvariant} to Gromov-Witten
theory.
Usually, the latter subject is studied for compact geometries.
For us, the open target consisting of $X^o$ relative to
the disjoint divisors $D_1^o,\ldots,D^o_n,D^o_{\out}$
is much more natural.
In order to define relative invariants for our open
geometries, properness of the associated
evaluations maps will be proven.
We will obtain another equivalence
\begin{equation}\label{vved}
N^{\text{hol}}_{\bf m}({\bf w}) =N^{\text{rel}}_{\bf m}({\bf w})
\end{equation}
where the term on the right is a relative Gromov-Witten count for
 $X^o$.

In fact, Gromov-Witten theory provides a direct interpretation of 
the {\em entire} summation of Theorem \ref{mtr}. We will
find $\log f_\fod$ is simply a summation of genus 0 Gromov-Witten
invariants of blow-ups of $X^o$ relative to $D^o_{\out}$. 
The result will be proven in Section \ref{ftvg}
from Theorem \ref{mtr} 
by the equivalences, 
the degeneration formula, and
exact  multiple cover calculations in Gromov-Witten theory.

\subsection{Desingularization}
\label{gws}
Following the notation of Section \ref{tcounts},
let $${\bf m}=(m_1,\ldots,m_n)$$ be an $n$-tuple of
  primitive vectors of $M$.
Let ${\bf w}=({\bf w}_1,\ldots, {\bf w}_n)$ be a $n$-tuple of weight vectors
$${\bf w}_i=(w_{i1},\ldots,w_{il_i})$$
 with
$$0\neq m_{\out}
=\sum_{i=1}^n |{\bf w}_i| m_i, \ \  m_{\out}=w_{\out}m'_{\out}$$
for $m'_{\out}\in M$ primitive.
We again treat the nondegenerate case where the ray
generated by $m_\out$ is distinct from the rays generated
by $-m_i$.
  
Let $\Sigma$ be the fan with rays generated by 
$$-m_1,\ \ldots,\ -m_n,\ m_{\out}\ .$$
Let $X$ be the toric surface over $\CC$ associated to
$\Sigma$ 
with toric divisors
$D_1,\ldots,D_n,D_{\out}$. Let
$\widetilde\Sigma$  be a refinement of $\Sigma$
satisfying the following two properties:
\begin{enumerate}
\item[(i)] the toric surface $\widetilde{X}$ associated
to $\widetilde{\Sigma}$ is nonsingular,
\item[(ii)] the  proper transforms of 
$D_1,\ldots,D_n,D_{\out}$
under the
birational toric  morphism 
$\widetilde X\rightarrow X$ are pairwise disjoint.
\end{enumerate}
Let 
$\widetilde D_1,\ldots,
\widetilde D_n,\widetilde D_{\out}\subseteq \widetilde{X}$ denote
the respective proper transforms of the divisors.

Let $\beta_{\bf w}\in H_2(\widetilde X,\ZZ)$ be as in \S
\ref{section-tropicaltoholomorphic}.
The class $\beta_{\bf w}$ is determined 
by intersection
numbers with toric divisors,
\begin{eqnarray*}
\widetilde D_i \cdot\beta_{\bf w} =|{\bf w}_i|,\quad&
\widetilde D_{\out}\cdot\beta_{\bf w}=w_{\out}
\end{eqnarray*}
and 
$D\cdot\beta_{\bf w}=0$ if 
$D\not\in\{\widetilde D_1,\ldots,\widetilde D_n,\widetilde D_{\out}\}$.
The full toric boundary  $B\subseteq \widetilde{X}$ 
represents the
anticanonical class of $\widetilde{X}$. Hence
$$\int_{\beta_{\bf w}} c_1(T_{\widetilde{X}}) = \sum_{i=1}^n
|{\bf w}_i|  + w_\out.$$

Let $X^o$ be the complement of the 0-dimensional torus orbits
in $X$, or equivalently, 
\[
X^o=\widetilde X\setminus \bigcup_D D,
\]
where the union is over all toric divisors $D$ of $\widetilde X$ with
$D\not\in \{\widetilde D_1,\ldots,\widetilde D_n,\widetilde D_{\out}\}$.
Every map
$$
\varphi:(\mathbb{P}^1,Q_{11},\ldots,Q_{1l_1},\ldots,
Q_{n1},\ldots,Q_{nl_n}, Q)\rightarrow X^o$$
in
$\foM(X^o,{\bf w})$
 yields a divisor $\varphi(\mathbb{P}^1)\cap D_i$ on each  $D_i$.

\begin{lemma} \label{nnmm}
The $n$ divisors 
\begin{equation*}
\varphi(\mathbb{P}^1)\cap D_1\subseteq D_1,\ \ldots,\ \varphi(\mathbb{P}^1)
\cap D_n\subseteq D_n
\end{equation*}
determine $\varphi(Q) \in D_\out$ up to $w_\out$ choices.
\end{lemma}

\proof
The full toric boundary 
$$\iota:B \rightarrow \widetilde{X}$$ is
a simple loop of $\PP^1$'s each meeting two others. Hence,
$B$ has arithmetic
genus 1 and
$\text{Pic}_0(B)\stackrel{\sim}{=}\CC^*$.
We have
$$\sum_{i=1}^n \varphi(\mathbb{P}^1)\cap D_i + w_\out \cdot \varphi(Q) 
= \iota^*( \beta_{\bf w} ) \ \in \text{Pic}(B),
$$ where $\beta_{\bf w}$ here is considered an element of 
$\text{Pic}(\widetilde{X})$.
We conclude $\varphi(Q)$ is a $w_\out$-root of 
$$\iota^*(\beta_{\bf w}) - \sum_{i=1}^n \varphi(\mathbb{P}^1)\cap D_i$$
in $\text{Pic}(B)$. There are exactly $w_\out$ possibilities for $\varphi(Q)$
in $D_\out$. 
\qed

\subsection{Stable relative maps}

We first describe a {\em partial} compactification 
by stable relative maps,
$$\foM(X^o,{\bf w}) \subseteq \overline{\foM}(X^o,{\bf w}).$$
equipped with an  evaluation morphism
\[
\ev:\overline{\foM}(X^o,{\bf w})\rightarrow\prod_{i=1}^n
(D_i^o)^{\length({\bf w}_i)} \, .
\]

Consider the geometry of $\widetilde{X}$ 
relative to the union of the divisors 
\begin{equation}\label{nttcc}
\widetilde D_1,\ldots,\widetilde D_n,\widetilde D_{\out} \subseteq \widetilde{X}.
\end{equation}
Since $\widetilde{X}$ is nonsingular and 
the nonsingular divisors \eqref{nttcc} are disjoint, the moduli space 
of stable relative maps to the geometry is a well-defined 
Deligne-Mumford stack \cite{IP,LR,L}.

Let
$\overline{\foM}(\widetilde X,{\bf w})$
denote the moduli space of
stable relative maps of genus 0 curves
 representing the class $\beta_{\bf{w}}$ 
with tangency conditions{\footnote{The relative
conditions specified here are simply the orders of the
tangencies at the marked points. 
The order of the tangency is the local intersection number. Order 1
is transverse intersection, order 2 is usual tangency, and so on.
The locations of the tangencies on the
relative divisors are not 
specified yet.
}} given by
\begin{enumerate}
\item[$\bullet$]
  the weight vector ${\bf{w}_i}$ along $\widetilde{D}_i$,
\item[$\bullet$] full tangency of order $w_\out$ at a 
single point along $\widetilde{D}_\out$.
\end{enumerate}
The moduli space $\overline{\foM}(\widetilde X,{\bf w})$ parameterizes
maps $\varphi$ from connected curves $C$ of arithmetic genus 0 with
at worst nodal singularities to 
 destabilizations{\footnote{A destabilization
along a relative divisor is obtained by attaching
a finite number of bubbles each of which is
a $\mathbb{P}^1$-bundle over the divisor.
We refer the reader to 
Section 1 of \cite{junex} for an introduction to the destabilizations
required for stable relative maps. Li uses the term
{\em expanded degeneration} for our destabilizations.}}
 $\pi$ of $\widetilde{X}$
along the relative divisors,
 \begin{equation}\label{bartt}
(C,Q_{11},\ldots,Q_{1l_1},\ldots,
Q_{n1},\ldots,Q_{nl_n}, Q)
\stackrel{\varphi}{\rightarrow} \widetilde{\foX} 
\stackrel{\pi}{\rightarrow} \widetilde{X} \ .
\end{equation}
The relative conditions require
$\varphi$ to be tangent to $\widetilde D_{\out}$ at $Q$ with order $w_{\out}$
and tangent to $\widetilde{D}_i$ at $Q_{ij}$ with order $w_{ij}$.

Let
$
\overline{\foM}(X^o,{\bf w})
$
be the open Deligne-Mumford substack of genus 0
stable relative maps to
$X^o\subseteq \tilde X$. More precisely,
$
\overline{\foM}(X^o,{\bf w})
$
consists of maps \eqref{bartt} for which the composition
$\pi \circ \varphi$ has 
image contained
in $X^o$.
Evaluation maps
\begin{equation}\label{vvqa}
\ev:\overline{\foM}(X^o,{\bf w})
\rightarrow \prod_{i=1}^n (D_i^o)^{l_i}.
\end{equation}
are naturally obtained.
The moduli space $\foM(X^o,{\bf w})$ defined in
Section \ref{hcounts}
is an open substack of $\overline{\foM}(X^o,{\bf w})$.

While $\overline{\foM}(X^o,{\bf w})$ is typically not
a proper stack, the evaluation maps \eqref{vvqa} 
are proper morphisms.

\begin{proposition} The evaluation map $\ev$ is proper. \label{pprr1}
\end{proposition}

\proof We use the valuative criterion for properness, so
let $R$ be a valuation ring with residue field $K$, and suppose
we are given the left-hand square in the following diagram:
\[
\begin{matrix}
T=\Spec K&\mapright{}&\overline{\foM}(X^o,{\bf w})&
\hookrightarrow&\overline{\foM}(\widetilde X,{\bf w})\\
\mapdown{}&&\mapdown{\ev}&&\mapdown{\widetilde\ev}\\
S=\Spec R&\mapright{}&\prod_{i=1}^n (D_i^o)^{l_i}&
\hookrightarrow &\prod_{i=1}^n \widetilde{D}_i^{l_i}
\end{matrix}
\]
where $\widetilde\ev$ is the corresponding evaluation map.
Since $\overline{\foM}
(\widetilde X,{\bf w})$ is proper,
$\widetilde\ev$ is certainly proper.
 Thus,
we obtain a unique family of stable relative 
maps
\[
\begin{matrix}
\shC&\mapright{\varphi}&\widetilde\foX_S&\mapright{\pi_S}&\widetilde X\times S\\
\mapdown{}&&\mapdown{}&&\\
S&\mapright{=}&S&&
\end{matrix}
\]
where $\widetilde\foX_S$ is a destabilization of $\widetilde{X} \times S$
 relative to the divisor
$\sum_{i=1}^n\widetilde D_i +\widetilde D_{\out}$.
We will show $\varphi$ is in fact a family of stable relative  maps
to $X^o$.

Let $0\in S$ be the closed point, and consider the 
morphism 
$$\pi_0\circ \varphi_0:\shC_0\rightarrow \widetilde X\ .$$
The marked points of $\shC_0$ all map to $ X^o$.
Suppose $\pi_0\varphi_0(\shC_0)$ intersects a toric
divisor $D\subseteq \widetilde{X}$ at a point of
$D\setminus X^o$.
Since the intersection number of $\shC_0$ with $D$ is 
accounted for in $X^o$, 
there {\em must} be an
irreducible component of $\shC_0$ dominating $D$.

If $\pi_0\circ
\varphi_0$ does not factor through $X^o\hookrightarrow \widetilde X$,
we have shown that
there must be an irreducible component $\widehat C$ of $\shC_0$ dominating
a toric divisor $D$ of $\widetilde X$.
 Let 
\[
\shC\mapright{\xi} \shC'\mapright{\psi} \widetilde{X}\times S
\]
be the Stein factorization of $\pi_S\circ\varphi$, so that
$\xi$ has connected fibers and $\psi$ is finite. Let
$
C'=\xi(\widehat C)$.
Since $\widehat C$ dominated $D$, $\xi$ cannot contract 
$\widehat C$. Hence $C'$
is an irreducible component of $\shC'_0$.

Let $E_1$ and $E_2$ be the two distinct
toric divisors
of $\widetilde X$ intersecting $D$ only at two distinct
 toric fixed points of $D$.
We will now show that there are irreducible components 
$$C_1', C_2' \subseteq \shC'_0$$
 intersecting $C'$ and dominating $E_1$ and $E_2$ respectively.

For $i\in \{ 1,2\}$,
let $\shL_i$ be the line bundle on $\widetilde X\times S$
corresponding to the divisor $E_i\times S$.
Let $s_i\in \Gamma(\shL_i)$
be a global section vanishing on $E_i\times S$. 
If $\psi^*(s_i)$
vanishes on an irreducible component $\shC'_i$ of $\shC'$,
then the generic fiber $(\shC'_i)_T$ of $\shC'_i\rightarrow S$ dominates 
$E_i\times T$, contradicting the 
assumption that the generic fiber $\shC_T$ of $\shC\rightarrow S$ maps
to $X^o$. 
Hence, $\psi^*(s_i)$ must define an effective
Cartier divisor $E_i'$ on $\shC'$. The support of
$\psi^*(s_i)$  is
codimension 1 in  $\shC'$ and intersects $C'$.

By assumption, all points of $\psi^{-1}(E_i\times
T)$ must be the images under $\xi$ of a marked point of $\shC_T$. Hence,
the closure of $\psi^{-1}(E_i\times T)$ in $\shC'$ consists of images
of marked points. Since the evaluation map
$S\rightarrow \prod D_i^{l_i}$ factors through $\prod (D_i^o)^{l_i}$
by assumption,{\footnote{The assumption covers all
the marked points except $Q$. However, the image
of $Q$ is constrained by the images of the other marked points by Lemma
\ref{nnmm} and therefore cannot intersect $C'$ even in the closure.}}
the closure of $\psi^{-1}(E_i\times T)$
in $\shC'$ cannot intersect $C'$. 
Thus $E_i'$ must contain
an irreducible component $C_i'$ of $\shC_0'$ intersecting $C'$. 
Then, $\psi(C_i')\subseteq E_i$.
Since $\psi$ is finite, $C_i'$ must dominate $E_i$.

By applying the 
above argument  repeatedly, replacing $C'$ with $C'_1$ or $C'_2$,
we find that either $\shC'_0$ has an infinite number of components
{\em or} $\shC'_0$ contains a cycle of components
dominating the union of toric divisors of $\widetilde X$.
The first possibility is absurd. In the  second case,
 the
genus of $\shC'_0$ would have to be at  least 1. Since $\xi|_{\shC_0}:\shC_0
\rightarrow \shC_0'$ has connected fibers, $\shC_0$ would also have
genus at least 1, contradicting the genus 0 assumption.

We have shown by contradiction that 
 $\varphi:\shC\rightarrow\widetilde{\foX}_S$ must in
fact be a family of stable relative maps to $X^o$.
\qed

\vspace{10pt}

Maps parameterized by  $\overline{\foM}(X^o,{\bf w})$
have a rather simple structure. Let 
 \begin{equation}\label{barttx}
(C,Q_{11},\ldots,Q_{1l_1},\ldots,
Q_{n1},\ldots,Q_{nl_n}, Q)
\stackrel{\varphi}{\rightarrow} {\foX}^o 
\stackrel{\pi}{\rightarrow} {X}^o \ 
\end{equation}
be an element of 
$\overline{\foM}(X^o,{\bf w})$.
Then, the
irreducible components 
$P\subseteq C$ are of exactly two distinct types
\begin{enumerate}
\item[(i)] $\pi\circ \varphi|_P:P \rightarrow X^o$ is a constant map, or
\item[(ii)] $\pi\circ \varphi|_P:P \rightarrow X^o$ is
an element of $\foM(X^o, {\bf w}^P)$ for
possibly different weight vectors ${\bf w}^P = ({\bf w}^P_1, \ldots
, {\bf w}^P_n)$, with not all
${\bf w}^P_i=\emptyset$.
 \end{enumerate}
The possibilities are restricted to (i) and (ii) because
no components of $C$ may both map to a destabilizing component
{\em and} surject onto the corresponding relative divisor.

The following result shows the domain $C$ is irreducible
if the locations of the relative conditions are specified
generically.

\begin{proposition}\label{pprr2}
The general fiber of 
\[
\ev:\overline{\foM}(X^o,{\bf w})
\rightarrow \prod_{i=1}^n (D_i^o)^{\length({\bf w}_i)}
\]
is contained in $\foM(X^o,{\bf w})$.
\end{proposition}

\proof
Let 
$C
\stackrel{\varphi}{\rightarrow} {\foX}^o 
\stackrel{\pi}{\rightarrow} {X}^o$
be an element of the general fiber of the evaluation map.
Consider first the component $P\subseteq C$ containing the point $Q$.
Suppose $P$ is of type (i). Let
$$B_1,B_2, \ldots,B_r \subseteq {\foX}^o$$
be the destabilizing components over $D^o_\out$ with
$$X^o \cap B_1 \neq \emptyset$$
and
$$B_i\cap B_{i+1}\not=\emptyset, \quad 1\le i\le r-1.$$
Each $B_i$ is a $\mathbb{P}^1$ bundle over $D^o_\out$
relative to two disjoint sections.
Since $P$ contains $Q$, the image of $P$ lies in $B_r$.
By the definition of a stable relative map, we conclude
\begin{enumerate}
\item[(a)]
The intersection of the image of $C$ with $\cup_{k=1}^r B_k$
is a single chain of $\mathbb{P}^1$-fibers which intersects
$X^o$ at a point $x\in D^o_\out$.
\item[(b)] There are at least two type (ii) components of $C$
intersecting $D^o_\out$ at $x$.
\end{enumerate}
The conclusion  (b) follows from the stability assumption
on the maps of components of $C$ to $\cup_{k=1}^r B_k$.

Let $P_1,P_2 \subseteq C$ be two components of type (ii) 
intersecting $D^o_\out$ at $x$.
The restrictions of $\pi \circ \varphi$ to $P_1$ and $P_2$
 determines
elements of 
$\foM(X^o, {\bf w}^{P_1})$ and 
$\foM(X^o, {\bf w}^{P_2})$ respectively.
By the genericity assumption, 
both restrictions lie in the general fibers of the  evaluation
maps 
\begin{eqnarray*}
\ev:\overline{\foM}(X^o,{\bf w}^{P_1})
& \rightarrow& \prod_{i=1}^n (D_i^o)^{\length({\bf w}^{P_1}_i)}, \\
\ev:\overline{\foM}(X^o,{\bf w}^{P_2})
&\rightarrow &\prod_{i=1}^n (D_i^o)^{\length({\bf w}^{P_2}_i)}.
\end{eqnarray*}
But the incidence point $x$ for the two restrictions is easily
seen {\em not} to agree for elements of general fibers
of the evaluation maps by Lemma \ref{nnmm}.


We have shown the component $P$ containing $Q$ must be of type (ii)
and no destabilization occurs over $D^o_\out$.
Any other component $P'\subseteq C$ of type (ii) cannot intersect $D^o_\out$
since all intersection have been accounted for.
But by Proposition \ref{Mdimension2} and the genericity statement,
no such curves $P'$ exist.
Hence, $P$ is the unique type (ii) component.

There are now only two possibilities for type (i) components
$\widehat{P}\subseteq C$.
First, $\widehat{P}$ could be contracted to a point of
$X^o$ away from the toric divisors. However, then $\widehat{P}$
would have to intersect $P$ in an least three points, violating
the genus condition on $C$.
Second, $\widehat{P}$ could be contracted to a point
of $D_i^o$. 
Then, let
$$B_1^i,B^i_2, \ldots,B^i_s \subseteq {\foX}^o$$
be the destabilizing components over $D^o_i$.
The component $\widehat{P}$ distinguishes a chain of
$\mathbb{P}^1$-fibers of
$\cup_{k=1}^s B^i_k$
which intersect
$X^o$ at a point $y\in D^o_i$.
Let $\widehat{P}_0\subseteq C$ be a component which
maps to the last element of the chain in $B^i_s$.
By the genericity of the fiber of the evaluation map,
$\widehat{P}_0$ is unique. The components of $C$ mapping to
the chain connect $\widehat{P}_0$  with $P$. These components
must themselves be a chain by the genus 0 condition on $C$.

The last argument applies to {\em every} component
$\widehat{P}$ of type (i) mapping to $D^o_i$.
Hence, we conclude the stability of relative map is violated 
over $D_i^o$.
Therefore, $C$ has a single type (ii) component and no
type (i) components. 
\qed

\subsection{Gromov-Witten theory for $X^o$}
The moduli space $\overline{\foM}(X^o,{\bf w})$
carries a perfect obstruction theory and a virtual class
$$[\overline{\foM}(X^o,{\bf w})]^{vir} \in A_*(
\overline{\foM}(X^o,{\bf w}), \mathbb{Q})$$
of dimension 
$$\int_{\beta({\bf w})} c_1(T_{X^o}) -1 - \sum_{i=1}^n \sum_{j=1}^{l_i} (w_{ij}-1) - (w_\out -1)
=\sum_{i=1}^n \text{length}({\bf w}_i),$$
see \cite{IP,LR,L}.

Since $\overline{\foM}(X^o,{\bf w})$ is not a
proper Deligne-Mumford stack, obtaining numerical
invariants from the virtual class is not always possible.
However, the properness of the evaluation map
\begin{equation}\label{vvqap}
\ev:\overline{\foM}(X^o,{\bf w})
\rightarrow \prod_{i=1}^n (D_i^o)^{l_i}
\end{equation}
proven in Proposition \ref{pprr1}
may be used to define the invariants we will need.

Select a point of $\prod_{i=1}^n (D_i^o)^{l_i}$, and let
$$\gamma:\text{Spec}(\CC) \rightarrow  \prod_{i=1}^n (D_i^o)^{l_i}$$
be the associated inclusion map.
The properness of \eqref{vvqap} implies
 $$F(\gamma) = \text{Spec}(\CC)  \times_ {\prod_{i=1}^n (D_i^o)^{l_i}}
\overline{\foM}(X^o,{\bf w})$$ is a proper Deligne-Mumford stack.
Define the relative Gromov-Witten invariant
\begin{equation}\label{fbttq}
N^{\text{rel}}_{\bf m}({\bf w})=
\int_{F(\gamma)}
\gamma^!( [\overline{\foM}(X^o,{\bf w})]^{vir})
\end{equation}
where integration is given as usual by proper push-forward to a point.

Since the codimension of the inclusion $\gamma$ equals
the virtual dimension of $\overline{\foM}(X^o,{\bf w})$,
definition \eqref{fbttq} is sensible.
The properness of the evaluation map \eqref{vvqap} implies
$N^{\text{rel}}_{\bf m}({\bf w})$ is {\em independent} of the 
point $\gamma$.

The invariant $N^{\text{rel}}_{\bf m}({\bf w})$ is the
virtual count of genus 0 maps to $X^o$ of class
$\beta_{\bf w}$ with tangency conditions
${\bf w}_i$ at {\em specified} points of
$D_i^o$ and
a tangency condition $w_\out$ at an unspecified point of $D_\out^o$.
 The actual moduli spaces of such maps may vary
with different choices of the locations of the tangency
conditions. For example, dimensions may jump. But the virtual
count $N^{\text{rel}}_{\bf m}({\bf w})$ does not change.

\begin{theorem} \label{pprr3}
$N^{\text{\em hol}}_{\bf m}({\bf w}) =N^{\text{\em rel}}_{\bf m}({\bf w}).$
\end{theorem}

\proof 
We can compute $N^{\text{rel}}_{\bf m}({\bf w})$ for a 
general point
$$\gamma:\text{Spec}(\CC) \rightarrow  \prod_{i=1}^n (D_i^o)^{l_i}.$$
By Proposition \ref{pprr2},
the fiber product $F(\gamma)$ embeds in 
${\foM}(X^o,{\bf w})$. By Proposition \ref{Mdimension},
${\foM}(X^o,{\bf w})$ is a nonsingular Deligne-Mumford stack 
of the expected dimension.
Hence, the virtual class
of ${\foM}(X^o,{\bf w})$ is the usual fundamental class.
By Bertini, the fiber product $F(\gamma)$ is a 
finite set of reduced points.
By definition, $N^{\text{rel}}_{\bf m}({\bf w})$ equals the 
 cardinality of $F(\gamma)$. On the other hand, 
the cardinality of $F(\gamma)$ is simply the degree of
\begin{equation*}
\ev:{\foM}(X^o,{\bf w})
\rightarrow \prod_{i=1}^n (D_i^o)^{l_i}
\end{equation*}
which is  $N_{\bf m}^{\text{hol}}({\bf w})$. \qed

\section{Formulas in the tropical vertex group}
\label{ftvg}
\subsection{Simple blow-ups} 
\label{sbu}

Let ${\bf P}=({\bf P}_1, \ldots, {\bf P}_n)$ be an $n$-tuple of
ordered partitions where 
$${\bf P}_i = p_{i1}+ \ldots +p_{i\ell_i}.$$
We continue to treat the nondegenerate case where the ray
generated by 
$$0\neq m_{\out}
=\sum_{i=1}^n |{\bf P}_i| m_i, \ \  m_{\out}=p_{\out}m'_{\out}$$
is distinct from the rays generated
by $-m_i$.
Following the notation of Section \ref{gws}, let $X^o$ be 
the corresponding open toric surface with divisors
$D^o_1, \ldots, D^o_n, D^o_\out$.

Select distinct points $x_{i1}, \ldots,x_{i\ell_i}\in D^o_i$
corresponding to the parts of ${\bf P}_i$.
Let
$$\nu: \widetilde{X}[{\bf P}] \rightarrow \widetilde{X}$$
be the blow-up
of $\widetilde{X}$ along all the points $x_{ij}$, and
let
$$X^o[{\bf P}] = \nu^{-1}(X^o).$$
Let $E_{ij} \subseteq \widetilde{X}[{\bf P}]$ be the exceptional divisor
over $x_{ij}$.
Let $\beta\in H_2(\widetilde{X},\mathbb{Z})$ be determined 
by intersection
numbers with toric divisors,
\begin{eqnarray*}
\widetilde D_i \cdot\beta =|{\bf P}_i|,\quad&
\widetilde D_{\out}\cdot\beta=p_{\out}
\end{eqnarray*}
and 
$D\cdot\beta=0$ if 
$D\not\in\{\widetilde D_1,\ldots,\widetilde D_n,\widetilde D_{\out}\}$.
Let $$\beta_{\bf P} = \nu^*(\beta) -\sum_{i=1}^{n} 
\sum_{j=1}^{\ell_i} p_{ij} [E_{ij}]\ \in H_2(\widetilde{X}
[{\bf P}], \mathbb{Z}) .$$

Let $\overline{\foM}(\widetilde{X}[{\bf P}]/\widetilde{D}_\out)$
be the moduli of 
stable relative maps of genus 0 curves
 representing the class $\beta_{\bf{P}}$
with full tangency of order $p_\out$ at an 
unspecified point of $\widetilde{D}_\out$.
Let
$$\overline{\foM}(X^o[{\bf P}]/{D}^o_\out) \subseteq
\overline{\foM}(\widetilde{X}[{\bf P}]/\widetilde{D}_\out)$$
be the open subspace of maps which avoid
$\widetilde{X}[{\bf P}] \setminus X^o[{\bf P}]$.
Unlike the theories considered in Section \ref{gwt},
the geometry 
$X^o[{\bf P}]/{D}^o_\out$ is relative to a {\em single} irreducible divisor.

\subsection{Gromov-Witten theory}
The virtual dimension of $\overline{\foM}(X^o[{\bf P}]/{D}^o_\out)$
is calculated by unwinding the definitions,

\begin{eqnarray*}
\text{virdim}_\CC \ 
\overline{\foM}(X^o[{\bf P}]/{D}^o_\out) & = & 
\int_{\beta_{{\bf P}}} c_1(T_{X^o[P]})  -1 - (p_\out -1) \\
& = & \sum_{i=1}^n |{\bf{P}}_i| -
\sum_{i=1}^{n} \sum_{j=1}^{\ell_i} p_{ij} + p_\out -1 - (p_\out -1)\\
& = & 0.
\end{eqnarray*}

\begin{proposition}\label{vrttq}
$\overline{\foM}(X^o[{\bf P}]/{D}^o_\out)$ is proper over $\text{\em Spec}(\CC)$. 
\end{proposition}

\proof The argument exactly follows the proof of Proposition \ref{pprr1}.
Let 
$$(C,Q) \stackrel{\varphi}{\rightarrow}  \foX^o[{\bf P}] \stackrel{\pi}
{\rightarrow} X^o[{\bf P}]$$
be an element of $\overline{\foM}(X^o[{\bf P}]/{D}^o_\out)$.
Here, $\foX^o[{\bf P}]$ is a destabilization of $X^o[{\bf P}]$ along
the relative divisor $D^o_\out$.
If we consider the composition
$$ \nu \circ \pi \circ \varphi : C \rightarrow X^o \subseteq \widetilde{X},$$
then the intersections with all the toric divisors of $\widetilde{X}$ are
specified --- occurring at the points $x_{ij}$ in $D_i^o$ and
the image of $Q$ in $D_\out^o$.
The limits in
$\overline{\foM}(\widetilde{X}[{\bf P}]/\widetilde{D}_\out)$
of 1-parameter families of elements of
$\overline{\foM}(X^o[{\bf P}]/{D}^o_\out)$
cannot meet points of 
$$\widetilde{X}[{\bf P}] \setminus X^o[{\bf P}].$$
Otherwise, an entire  strict transform of
a toric divisor of $\widetilde{X}$ would lie in the image of the
limit, and the loop construction of the proof of Proposition
\ref{pprr1} can be made.
\qed
\vspace{10pt}

We define the Gromov-Witten invariant $N_{\bf m}[{\bf P}]$
by the usual integral,
$$N_{\bf m}[{\bf P}] 
= \int_{[ \overline{\foM}(X^o[{\bf P}]/{D}^o_\out)]^{vir}}
1.$$
The properness Proposition
\ref{vrttq}  holds relatively as the centers $x_{ij}\in D_i^o$ of the
blow-ups of $\widetilde{X}$ are moved along $D_i^o$. 
Hence,
the integral $N_{\bf m}[{\bf P}]$ is {\em independent} of the
locations of the distinct points $x_{ij}\in D_i^o$.

\subsection{Degeneration} \label{ff677}
Let $\CC$ be the affine line.
Let 
\begin{equation}\label{bbqz}
X^o \times \CC \rightarrow \CC
\end{equation}
be the trivial family over $\CC$.
The standard degeneration of $X^o$ 
relative to $D_1^o, \ldots, D_n^o, D^o_\out$ over $\CC$,
\begin{equation}\label{udeg}
\epsilon: \mathcal{F} \rightarrow \CC,
\end{equation}
is obtained by blowing-up 
the loci $D_1^o, \ldots D_n^o, D^o_\out$ over $0\in \CC$
in \eqref{bbqz}.
After blowing-up the sections of $\epsilon$
corresponding  to the points $x_{ij}$, we 
obtain
\begin{equation} \label{vvvv}
\epsilon_{\bf P}: \mathcal{F}[{\bf P}] \rightarrow \CC.
\end{equation}
For $\xi\neq 0$,
we have $\epsilon_{\bf P}^{-1}(\xi)\stackrel{\sim}{=} X^o[{\bf P}]$. The
special fiber has a different form
$$\epsilon_{\bf P}^{-1}(0)= X^o \cup \bigcup_{i=1}^n B^o[{\bf P}_i],$$
where $B^o[{\bf P}_i]$ is the blow-up of
$\mathbb{P}(\mathcal{O}_{D_i^o} \oplus \text{Norm}_{X^o/D_i^o})$
at the limits 
$$\overline{x}_{i1}, \ldots, \overline{x}_{i\ell_i}$$
of the points $x_{ij}$.

The moduli spaces of stable relative maps to
the fibers of $\epsilon_{\bf P }$ over $\xi \neq 0$ are  canonically
isomorphic to 
$\overline{\foM}(X^o[{\bf P}]/{D}^o_\out)$.
The limits of stable relative maps to 
$\epsilon_{\bf P }^{-1}(\xi)$
as $\xi \rightarrow 0$ are well-understood
in relative Gromov-Witten theory \cite{IP,LR,L}. 
The limit is a union of stable maps to the relative geometries
determined by the components of the special fiber over 0.
A priori, the limit may {\em leave} the open subspaces we
are considering. However, the properness argument (used twice
already) easily shows the 
limit is a union of
stable maps to the following $n+1$ open relative geometries 
$$X^o/ D_1^o \cup \ldots D_n^o \cup D_\out^o\ , \ \ 
B^o[{\bf P}_1]/D_1^o, \ \ \ldots, \ \ B^o[{\bf P}_n]/D_n^o, $$
with matching conditions along the common relative divisors.

Relative maps to $B^o[{\bf P}_i]/D_i^o$ are simple to describe.
Let 
$\overline{E}_{ij} \subseteq B^o[{\bf P}_i]$
be the exceptional divisor of 
$$B^o[{\bf P}_i]  {\rightarrow} 
\mathbb{P}(\mathcal{O}_{D_i^o} \oplus \text{Norm}_{X^o/D_i^o})$$
over $\overline{x}_{ij}$. 
Let 
$C_{ij} \subseteq B^o[{\bf P}_i]$
be the strict transform to $B^o[{\bf P}_i]$ of the unique
fiber of 
$$\mathbb{P}(\mathcal{O}_{D_i^o} \oplus \text{Norm}_{X^o/D_i^o}) \rightarrow D_i^o$$
 containing $\overline{x}_{ij}$. 
A relative map to $B^o[{\bf P}_i]/D_i^o$
with a connected domain and  class $d[C_{ij}]$ 
must be a $d$-fold multiple cover of $C_{ij}$.

We can calculate all the invariants of $B^o[{\bf P}_i]$ relative 
to $D_i^o$ which we will require.
Let
$$\overline{\foM}( B^o[{\bf P}_i]/ D_i^o, 
d)_j$$
be the moduli space of genus 0 relative maps of class $d[C_{ij}]$
with full tangency of order $d$ along $D_i^0$.
Let $\overline{\foM}(C_{ij}/\infty, d)$
be the moduli space of genus 0 relative maps of degree $d$
with ramification of order $d$ over $\infty =C_{ij}\cap D_i^o$.
As
spaces, 
\begin{equation} \label{ww34}
\overline{\foM}( B^o[{\bf P}_i]/ D_i^o, 
d)_j =
\overline{\foM}(C_{ij}/\infty, d)\  .
\end{equation}
However, the two moduli spaces  \eqref{ww34} carry 
obstruction theories which differ by
$$H^1(C, \varphi^*(\text{Norm}_{B^o[{\bf{P}}_i]/C_{ij}})),$$
at the moduli point  
$[\varphi:C\rightarrow C_{ij}]$.
The virtual dimension of 
$\overline{\foM}( B^o[{\bf P}_i]/ D_i^o, 
d)_j$
is easily seen to be 0. Let
$$R_d  =  \int_{[\overline{\foM}( B^o[{\bf P}_i]/ D_i^o, 
d)_j]^{vir}} 1 .$$

\begin{proposition} \label{rddd}
$R_d  = \frac{(-1)^{d-1}}{d^2}$.
\end{proposition}

\proof
The integral $R_d$ is 
\begin{eqnarray*}
R_d
& = & \int_{[\overline{\foM}(C_{ij}/\infty, d)]^{vir}}
e( H^1(C, \varphi^*(\text{Norm}_{B^o[{\bf{P}}_i]/C_{ij}})))\\
& = & \int_{[\overline{\foM}({\mathbb P}^1/\infty, d)]^{vir}}
e( H^1(C, \varphi^*({\mathcal O}_{{\mathbb P}^1}(-1))))\\
& = & \frac{(-1)^{d-1}} {d^2}\ .
  \\
\end{eqnarray*}
In the first and second lines, $e(V)$ denotes the Euler class
(top Chern class) of the vector bundle on the
moduli space of stable relative maps with fibers defined
by $V$.
The evaluation in the last equality is the 
genus 0 
part{\footnote{
The local Gromov-Witten theory of curves is proven in \cite{Brypan2}
to determine a TQFT. As such, the entire theory is
specified by values on 
the cap, the tube, and the pair of pants.
The integral $R_d$ arises in the
genus 0 part of the cap contribution. In fact, the full
genus $g$ cap is computed in \cite{Brypan} by Hodge integral
techniques.}}
of Theorem 5.1 of  \cite{Brypan}. \qed
{\vspace {10pt}}

An $n$-tuple 
  ${\bf w}$ of weight vectors is of the same 
{\em type} as 
an $n$-tuple ${\bf P}$ of ordered partitions
if
$|{\bf w}_i|= |{\bf P}_i|$
for all $i$.
A set partition of 
${\bf w}_i$ {\em compatible} with ${\bf P}_i$ is
a disjoint union{\footnote{Recall the length of ${\bf{w}}_i$ is $l_i$ and
the length of ${{\bf P}_i}$ is $\ell_i$. The set
$I_j$ is allowed to be empty if $p_{ij}=0$.}}
$$I_1 \cup  \ldots  \cup I_{\ell_i} = \{ 1, \ldots, l_i\}$$
satisfying
$$p_{ij} = \sum_{r\in I_j} w_{ir}$$
for all $j$.
Let 
$$R_{{\bf P}_i|{\bf w}_i} = \sum_{I_\bullet}  \prod_{j=1}^{l_i} R_{w_{ij}}$$
where the sum is over all 
set partitions $I_\bullet$ of 
${\bf w}_i$ compatible with ${\bf P}_i$.

The degeneration formula for relative Gromov-Witten
theory  applied to our setting yields the following
result.

\begin{proposition} \label{yy7}
We have 
$$N_{\bf m}[{\bf P}] =
\sum_{\bf w} 
N_{\bf m}^{\text{\em rel}}({\bf w})  \cdot 
\prod_{i=1}^n  
\frac{\prod_{j=1}^{l_i} w_{ij}}
{|\text{\em Aut}({\bf w}_i)|} \cdot
R_{{\bf P}_i|{\bf w}_i}
 $$
where the sum is over all $n$-tuples of weight vectors
${\bf w}$ of the same type as ${\bf P}$.
\end{proposition}

\proof
We simply apply the degeneration formula of relative 
Gromov-Witten theory \cite{IP,LR,L} to the family $\epsilon_{\bf P}$
of \eqref{vvvv}. Our properness results show the
formula both makes sense and is valid in the open geometry
at hand.
The left term $N^{\text{rel}}_{\bf m}({\bf w})$
is the contribution of $X^o$.
The numerator
$$\prod_{i=1}^n \prod_{j=1}^{l_i} w_{ij}$$
is the standard multiplicity in the degeneration
formula.
The ratio
$$ \frac{R_{{\bf P}_i|{\bf w}_i}}
{|\text{Aut}({\bf w}_i)|}$$
is exactly the correct automorphism weighted
contribution of the component 
$B^o[{\bf P}_i]$. 
\qed

\subsection{First commutator formulas}
\label{fcf}
We can now prove our first formulas in the tropical vertex group.
Let $\fod_1 \neq \fod_2$ be two lines
through the origin in $M_{\mathbb{R}}$.
 Let ${\bf m}=(m_1,m_2)$ be the two associated primitive 
vectors. Let
\begin{equation}\label{zzxbv}
f_{\fod_1}= 
\prod_{\xi=1}^{\ell_1} (1+ s_{\xi} z^{m_1}) , \ \ f_{\fod_2} 
= \prod_{\xi=1}^{\ell_2}(1+t_{\xi} z^{m_2})
\end{equation}
be functions over the
complete ring $\CC[[s_1, \ldots, s_{\ell_1},t_1,\ldots, t_{\ell_2}]]$.
Consider the scattering diagram 
\begin{equation}\label{gg234}
\foD=\{(\fod_1,f_{\fod_1}), (\fod_2,f_{\fod_2})\}\ .
\end{equation}

Let $\Scatter(\foD)$ be the unique minimal scattering diagram
 (obtained by adding
rays to $\foD$)  for which the path ordered product around the 
origin is trivial.
Let $m'_\out\in M$ be a primitive vector. What is the
associated function $f_{m'_\out}$?

We may write the scattering diagram \eqref{gg234} equivalently with $\ell_1+\ell_2$
lines,
\begin{equation*}
\foD=\{(\fod_1, 1+s_1 z^{m_1} ),\ldots,(\fod_1, 1+s_{\ell_1} z^{m_1} ),
(\fod_2, 1+t_1 z^{m_2}),\ldots,
(\fod_2, 1+t_{\ell_2} z^{m_2} )\} \ .
\end{equation*}
Then, $\foD$ is standard by Definition \ref{ssss}, and 
Theorem \ref{mtr} determines $f_{m'_\out}$ in terms of tropical geometry.
There are three possibilities
\begin{enumerate}
\item[(i)] $m'_\out \notin \ \mathbb{Q}_{\geq 0} m_1 + \mathbb{Q}_{\geq 0} m_2$,
\item[(ii)] $m'_\out =m_1$ or  $m'_\out= m_2$,
\item[(iii)] $m'_\out \in \ \mathbb{Q}_{> 0} m_1 + \mathbb{Q}_{>0} m_2$,
\end{enumerate}
The outcomes in cases (i) and (ii) are straightforward. By Theorem \ref{mtr}, 
 $f_{m'_\out}=1$  in case (i) and $f_{m'_\out}$ equals either
$f_{\fod_1}$ or $f_{\fod_2}$ respectively. The interesting case is (iii).

For ordered partitions ${\bf P}_1= p_{11}+ \ldots + p_{1\ell_1}$
and
${\bf P}_2= p_{21}+ \ldots + p_{2\ell_2}$,
let
$$s^{{\bf P}_1} = \prod_{\xi=1}^{\ell_1} s_\xi^{p_{1\xi}}, \ \ \
t^{{\bf P}_2} = \prod_{\xi=1}^{\ell_2} t_\xi^{p_{2\xi}}$$
be the corresponding monomials.

\begin{theorem}\label{hhhp} If $m'_\out \in \ \mathbb{Q}_{> 0} m_1 + \mathbb{Q}_{>0} m_2$, then 
$$\log f_{m'_\out} = \sum_{k=1}^\infty 
\ \sum_{\bf{P}=({\bf P}_1, {\bf P}_2)}
        k\ N_{\bf{m}}[{\bf{P}}]\ s^{{\bf P}_1}\ t^{{\bf P}_2}  \
z^{km'_\out}$$
where the sum is over all ordered partitions ${\bf P}_i$
of length $\ell_i$ satisfying
$$|{\bf P}_1|m_1 + |{\bf P}_2| m_2 = k m'_\out.$$
\end{theorem}

\proof
To apply Theorem \ref{mtr}, we first calculate the Taylor series
$$\log (1 + s_\xi z^{m_1}) = \sum_{d=1}^\infty  d\ \frac{(-1)^{d-1}}{d^2} z^{dm_1} s_\xi^d,$$
$$\log (1 + t_\xi z^{m_2}) = \sum_{d=1}^\infty  d\  \frac{(-1)^{d-1}}{d^2} z^{dm_2} t_\xi^d.$$
In  both of the above cases, 
the $a$ coefficients of Theorem \ref{mtr} match the
relative invariants $R_d$ computed in Proposition \ref{rddd},
$$a_{\xi dd}=R_d = \frac{(-1)^{d-1}}{d^2}.$$
By Theorems \ref{mtr}, \ref{N=Ntrop}, and \ref{pprr3},
$$\log f_{m'_\out} = \sum_{k=1}^\infty \
\sum_{\bf{P}=({\bf P}_1, {\bf P}_2)}
\sum_{\bf{w}}
        k\ N_{\bf{m}}^{\text{rel}}({\bf{w}})\ 
\left( \prod_{i=1}^2
\frac{\prod_{j=1}^{l_i} w_{ij}}
{|\text{Aut}({\bf w}_i)|} \cdot
R_{{\bf P}_i|{\bf w}_i} \right)
\
s^{{\bf P}_1}\ t^{{\bf P}_2}  \
z^{km'_\out} \ ,
$$
where the second sum is over all
 ordered partitions ${\bf P}=({\bf P}_1, {\bf P}_2)$
of lengths $(\ell_1, \ell_2)$ satisfying
$$|{\bf P}_1|m_1 + |{\bf P}_2| m_2 = k m'_\out$$
and the third sum is over all 
weight vectors ${\bf w}=({\bf w}_1, {\bf w}_2)$
of the same type as ${\bf{P}}$.
The result then follows from Proposition \ref{yy7}.
\qed
\vspace{10pt}

Theorem \ref{hhhp} computes the commutators of
$\theta_{\fod_1,f_{\fod_1}}$ and $\theta_{\fod_2,f_{\fod_2}}$
via slope ordered products in the tropical vertex group
in the form discussed in Section \ref{ttss}. 
Theorem \ref{hhhp} specializes to
Theorem \ref{fff}, but
is a much better statement.  By working over
the ring $\CC[[s_1, \ldots, s_{\ell_1},t_1,\ldots, t_{\ell_2}   ]]$, 
we see computing the
commutator is precisely {\em equivalent} to calculating {all} 
the Gromov-Witten invariants $N_{\bf m}[({\bf P}_1, {\bf P}_2)]$
for ordered partitions ${\bf P}_i$ of length $\ell_i$.

Of course, functions $f_1$ and $f_2$ may not be in 
 the form specified by \eqref{zzxbv}. The
general commutator formula is expressed in terms
of orbifold blow-ups of toric surfaces.

\subsection{Orbifold blow-ups}
Let $D\subseteq S$ be a nonsingular divisor contained
in a surface. Let $p\in D$.
Let $p_{D}^r$ be the unique length $r$ subscheme of
$D$ with support $p$. Let 
$$S_r \rightarrow S$$
be the blow-up of $S$ along $p_{D}^r$.
Of course, $S_1$ is the usual blow-up of $S$ along $p$.
For $r\geq 2$, $S_r$ has a unique 
$A_{r-1}$-singularity. Hence, $S_r$ admits a
unique structure as a  nonsingular Deligne-Mumford stack 
$\mathcal{S}_r \rightarrow S_r$.
We call the composition
$$\delta_r:\mathcal{S}_r\rightarrow S$$
the {\em $r$-orbifold blow-up} of $S$ along $(p,D)$

The exceptional divisor $E\subseteq {\mathcal {S}}_r$ of $\delta_r$
is a
$\mathbb{P}^1$ with a single $r$-fold stack point
lying above the original $A_{r-1}$-singularity.
The self-intersection is
$$[E]^2 = -\frac{1}{r}.$$

A {\em graded} partition consists of a finite sequence
$({\bf P}^1, {\bf P}^2, \ldots, {\bf P}^{d})$
of ordered partitions for which {\em every} part of
${\bf P}^r$ is divisible by  $r$. The length of a graded
partition is given by a $d$-tuple
$(\ell^1, \ldots, \ell^{d})$. The size of a graded
partition is
$\sum_{r=1}^d |{\bf P}^r|$.

Let ${\bf G}=({\bf G}_1, \ldots, {\bf G}_n)$
be an $n$-tuple of graded partitions,
 $${\bf G}_i= ({\bf P}^1_i, {\bf P}^2_i, \ldots, {\bf P}^{d_i}_i).$$
We treat the nondegenerate case where the ray
generated by 
$$0\neq m_{\out}
=\sum_{i=1}^n  |{\bf G}_i| m_i, 
\ \  m_{\out}=g_{\out}m'_{\out}$$
is distinct from the rays generated
by $-m_i$.
Let $X^o$ be 
the corresponding open toric surface with divisors
$D^o_1, \ldots, D^o_n, D^o_\out$.

Select a distinct point $x_{ij}^r$ of $D^o_i$
for each part of $p_{ij}^r$ of ${\bf P}_i^r$.
Let
$$\nu: \widetilde{X}[{\bf G}] \rightarrow \widetilde{X}$$
be obtained 
from $\widetilde{X}$ by taking 
$r$-orbifold blow-ups  along $(x^r_{ij}, D_i^o)$ for all $i,r,j$.
Let
$$X^o[{\bf G}] = \nu^{-1}(X^o).$$

Let $E^r_{ij} \subseteq \widetilde{X}[{\bf G}]$ be the exceptional divisor
over $x_{ij}^r$.
Let $\beta\in H_2(\widetilde{X},\mathbb{Z})$ be determined 
by intersection
numbers with toric divisors,
\begin{eqnarray*}
\widetilde D_i \cdot\beta =|{\bf G}_i|,\quad&
\widetilde D_{\out}\cdot\beta=g_{\out}
\end{eqnarray*}
and 
$D\cdot\beta=0$ if 
$D\not\in\{\widetilde D_1,\ldots,\widetilde D_n,\widetilde D_{\out}\}$.
Let $$\beta_{\bf G} = \nu^*(\beta) -\sum_{i=1}^{n}
\sum_{r=1}^{d_i}  
\sum_{j=1}^{\ell_i^r}  p^r_{ij} [E^r_{ij}]\ \in H_2(\widetilde{X}
[{\bf G}], \mathbb{Z}) .$$

Let $\overline{\foM}(\widetilde{X}[{\bf G}]/\widetilde{D}_\out)$
be the moduli of 
stable relative maps of genus 0 curves
 representing the class $\beta_{\bf{G}}$
with full tangency of order $g_\out$ at an 
unspecified point of $\widetilde{D}_\out$.
The moduli space of stable maps to orbifold targets is defined in
\cite{AGV,CR}.
Since the relative divisor $\widetilde{D}_\out$ does
not meet the orbifold points of $\widetilde{X}[{\bf G}]$,
there is no difficulty in defining the moduli space of 
relative maps.{\footnote{Orbifold stable maps are allowed also to
have nontrivial prescribed orbifold structure on
the domain. Our maps to $\widetilde{X}[{\bf G}]$ have no such prescribed
orbifold structure.}}
Let
$$\overline{\foM}(X^o[{\bf G}]/{D}^o_\out) \subseteq
\overline{\foM}(\widetilde{X}[{\bf G}]/\widetilde{D}_\out)$$
be the open subspace of maps which avoids
$\widetilde{X}[{\bf G}] \setminus X^o[{\bf G}]$.

Since the curve class $\beta_{\bf G}$ accounts for all
the intersections of the image of the
map in $\widetilde{X}$, the argument for properness is still
valid.

\begin{proposition}\label{vrttq5}
$\overline{\foM}(X^o[{\bf G}]/{D}^o_\out)$ is proper over $\text{\em Spec}(\CC)$. 
\end{proposition} 

The virtual dimension of $\overline{\foM}(X^o[{\bf G}]/{D}^o_\out)$
is 0.
Define the orbifold Gromov-Witten invariant as before
$$N_{\bf m}[{\bf G}] 
= \int_{[ \overline{\foM}(X^o[{\bf G}]/{D}^o_\out)]^{vir}}
1.$$
The integral $N_{\bf m}[{\bf G}]$ is {\em independent} of the
locations of the distinct points $x^r_{ij}\in D_i^o$.

\subsection{Full commutator formulas}
Let $\fod_1 \neq \fod_2$ be two lines
through the origin in $M_{\mathbb{R}}$. 
Let ${\bf m}=(m_1,m_2)$ be the two associated primitive 
vectors. Let
\begin{equation}\label{zzxb}
f_{\fod_1}= \prod_{r=1}^{d_1}
\prod_{\xi=1}^{\ell^r_1} (1+ s_{\xi}^r z^{r m_1}) , \ \ f_{\fod_2} 
= \prod_{r=1}^{d_2}
\prod_{\xi=1}^{\ell^r_2}(1+t_{\xi}^r z^{rm_2})
\end{equation}
over the complete ring  $\CC[[s_\bullet^\bullet,t_\bullet^\bullet]]$ 
in all the variables 
$s_{\xi}^r, t_{\xi}^r$.

For graded partitions ${\bf G}_1= ({\bf P}^1_1,\ldots, {\bf P}_1^{d_1})$
and ${\bf G}_2= ({\bf P}^1_2,\ldots, {\bf P}_2^{d_2})$,
let
$$s^{{\bf G}_1} = \prod_{r=1}^{d_1}
\prod_{\xi=1}^{\ell^r_1}\ (s^r_\xi)^\frac{p^r_{1\xi}}{r}, \ \ \
t^{{\bf G}_2} = \prod_{r=1}^{d_2}
\prod_{\xi=1}^{\ell^r_2} \ (t^r_\xi)^\frac{p^r_{2\xi}}{r}
$$
be the corresponding monomials.

Again, we consider the unique minimal scattering diagram
$\Scatter (\foD)$ associated to 
\begin{equation*}\label{gg2347}
\foD=\{(\fod_1,f_{\fod_1}), (\fod_2,f_{\fod_2})\}\ .
\end{equation*}

\begin{theorem}\label{hhhpx} 
If $m'_\out \in \ \mathbb{Q}_{> 0} m_1 + \mathbb{Q}_{>0} m_2$,
$$\log f_{m'_\out} = \sum_{k=1}^\infty 
\sum_{\bf{G}=({\bf G}_1, {\bf G}_2)}
        k\ N_{\bf{m}}[{\bf{G}}]\ s^{{\bf G}_1}\ t^{{\bf G}_2}  \
z^{km'_\out}$$
where the sum is over all graded partitions ${\bf G}_1$ 
of length
$(\ell_1^1, \ldots, \ell_1^{d_1})$ and 
${\bf G}_2$
of length
$(\ell_2^1, \ldots, \ell_2^{d_2})$
 satisfying
$$|{\bf G}_1|m_1 + |{\bf G}_2| m_2 = k m'_\out.$$
\end{theorem}

\proof
We follow the proof of Theorem \ref{hhhp}.
After factoring the diagram $\foD$ into 
$\sum_{r=1}^{d_1} \ell_1^r + \sum_{r=1}^{d_2} \ell_2^r$
lines, we will match
the formula for $\log f_{m'_\out}$ from Theorem \ref{mtr}
with a degeneration calculation of the
orbifold Gromov-Witten invariants of  
$X^o[({\bf{G}}_1,{\bf{G}}_2)]/D_\out^o$.

The geometric setting here is just as before. First, the
degeneration 
$$\epsilon: {\mathcal {F}} \rightarrow \CC$$
defined in \eqref{udeg} is taken.
The points
$x_{1\xi}^r \in D_1^o$ and $x_{2\xi}^r \in D_2^o$
specialize to points of 
\begin{equation}\label{hh12}
\mathbb{P}(\mathcal{O}_{D_1^o} \oplus \text{Norm}_{X^o/D_1^o}) \ \ 
{\text{ and }} \ \ 
\mathbb{P}(\mathcal{O}_{D_2^o} \oplus \text{Norm}_{X^o/D_2^o})
\end{equation}
respectively.
The projective bundles \eqref{hh12} each contain two distinguished
sections. The first section is the limit in $\mathcal{F}$ of
the divisor $D_i^o$ and
carries the limits of the $x_{i\xi}^r$.
The second section
meets
$X^o$.
There is no difficulty in taking the $r$-orbifold
blow-ups relative to $\epsilon$. The underlying coarse space is
obtained by blowing-up the families of canonically defined
nilpotent subschemes, and the stack structure is uniquely determined.  
The resulting
family,
$$\epsilon_{\bf G}: {\mathcal F}[{\bf G}] \rightarrow \CC,$$
has orbifold structure in the 
total space ${\mathcal F}[{\bf G}]$ 
with support disjoint from
the relative divisors in the special fiber $\epsilon_{\bf G}^{-1}(0)$.
Therefore, the usual degeneration formulas in relative
Gromov-Witten theory hold unchanged.

A crucial fact used in Theorem \ref{hhhp} is the surprising
match between the $a$ coefficients and the multiple
cover contributions $R_d$ of Proposition \ref{rddd}.
We need to be even luckier now. Since
$$\log (1 + s^r_\xi z^{rm_1}) 
= \sum_{d=1}^\infty  dr\ \frac{(-1)^{d-1}}{rd^2} z^{drm_1} (s^r_\xi)^d,$$
$$\log (1 + t^r_\xi z^{rm_2}) 
= \sum_{d=1}^\infty  dr\  \frac{(-1)^{d-1}}{rd^2} z^{drm_2} 
(t^r_\xi)^d,$$
the $a$ coefficients in both cases are
$$a^r_{\xi, d, rd} = \frac{(-1)^{d-1}}{rd^2}.$$
These exactly match the orbifold multiple cover
calculation in Proposition \ref{sddd} below.

The remaining steps are identical to those taken in
the proof of Theorem \ref{hhhp}. The formula of Theorem
\ref{mtr} matches exactly 
the degeneration computation. 
\qed

\vspace{10pt}

Let $P$ be the fiber of the projective bundle 
$$\mathbb{P}(\mathcal{O}_{D_i^o} \oplus \text{Norm}_{X^o/D_i^o})
\rightarrow D_i^o$$
containing the limit $\overline{x}_{i\xi}^r$.
Let $C_{i\xi}^r$ be the strict transform of $P$ after taking the
$r$-orbifold blow-up along $(\overline{x}_{i\xi}^r, D_i^o)$.
Then,
 $$C_{i\xi}^r \stackrel{\sim}{=}{\mathbb{P}}^1[r,1]$$
 with a single orbifold point
of order $r$ at $0$.
The normal bundle of $C_{i\xi}^r$ in the $r$-orbifold blow-up is
simply $\mathcal{O}_{C_{i\xi}^r}(-[0]/{\mathbb{Z}}_r)$ of degree $-\frac{1}{r}$.
For any map{\footnote{We mean here
a representable map. The domain $C$ has
no stack structure.}} from a genus 0 curve 
$$\varphi:C \rightarrow C_{i\xi}^r,$$
the degree of $\varphi$ must be a multiple of $r$.
Let
$$R^r_d  = 
 \int_{[\overline{\foM}( {\mathbb P}^1[r,1]/ \infty,rd)]^{vir}} 
e(H^1(C,\varphi^*(\mathcal{O}_{C_{i\xi}^r}(-[0]/{\mathbb{Z}_r})))) ,$$
where $\overline{\foM}( {\mathbb P}^1[r,1]/ \infty,rd)$
is the moduli space of genus 0 stable relative maps of degree $rd$ 
with full ramification $rd$
over $\infty$. 

\begin{proposition} \label{sddd}
$R^r_d  = \frac{(-1)^{d-1}}{rd^2}$.
\end{proposition}

\proof
There is a $\CC^*$ action on ${\mathbb P}^1[r,1]$ with tangent
weights $[\frac{1}{r},-1]$ at the fixed points $[0,\infty]$.
The orbifold line bundle 
$\mathcal{O}_{C_{i\xi}^r}(-[0]/{\mathbb{Z}_r})$
is canonically linearized with fiber weights $[-\frac{1}{r},0]$
over the respective fixed points.

We compute the integral $R^r_d$ via the induced $\CC^*$
action on $\overline{\foM}( {\mathbb P}^1[r,1]/ \infty,rd)$
and the Bott residue formula.
There are many $\CC^*$-fixed loci in the moduli space. However,
the above linearization leads to the vanishing{\footnote{See
the proof of Theorem 5.1 in \cite{Brypan} where the
same vanishing is explained.}}
 of {\em all}
contributions except for the single $\CC^*$-fixed Galois
cover
$$\varphi: C \rightarrow {\mathbb P}^1[r,1]$$
of degree $rd$.
The product of the weights of $\CC^*$ on 
$H^1(C,\varphi^*(\mathcal{O}_{C_{i\xi}^r}(-[0]/{\mathbb{Z}_r})))$
are
$$ \prod_{i=1}^{d-1} -\frac{ i}{rd} $$
The product of the weights of $\CC^*$ on the tangent space to
$[\varphi]$ in 
$\overline{\foM}( {\mathbb P}^1[r,1]/ \infty,rd)$ is
$$\frac{\prod_{i=1}^{d}\frac{i}{rd}}{\frac{1}{rd}}$$
where the bottom factor is obtained from reparameterization 
over $0\in {\mathbb{P}}^1[r,1]$. The above weight calculations, obtained
from the $\CC^*$-equivariant geometry of basic orbifold line bundle on
${\mathbb P}^1[r,1]$,
are standard. See
Section 2.2 of \cite{jpt} for a detailed treatment.

By the Bott residue formula,
$R^r_d$ equals the ratio of the above weights together with a stack
automorphism factor of $\frac{1}{rd}$, 
$$R^r_d = \frac{1}{rd} \frac{\prod_{i=1}^{d-1} -\frac{ i}{rd}}
{\frac{\prod_{i=1}^{d}\frac{i}{rd}}{\frac{1}{rd}}} = \frac{(-1)^{d-1}}{rd^2}.$$
Our calculation is just a minor modification of Theorem 5.1 of \cite{Brypan}.
\qed

{\vspace {10pt}}

Consider the tropical vertex group over the ring $\CC[[t]]$.
A function 
$$f= 1 + tz^m \cdot g(z^m,t)\ \ \ g \in \CC[z^m][[t]]$$ 
attached to the ray with primitive $m$ can always be factored
as
\begin{equation}\label{ffak}
f = (1+s_1^1 z^m)(1+s_1^2 z^{2m})(1+s_1^3 z^{3m}) \ldots 
\end{equation}
for $s_1^r\in t^{n_r}\CC[[t]]$ with 
$$\lim_{r\rightarrow \infty} n_r = \infty.$$
To any finite order in $t$, a suitable finite 
truncation of the factorization \eqref{ffak} suffices
for any calculation.
Therefore, we view Theorem \ref{hhhpx} as an ordered product formula for
an arbitrary commutator in the tropical vertex group.

We have written Theorem \ref{hhhpx} for functions \eqref{zzxb} in slightly
more complicated
form than \eqref{ffak} to capture all of the possible
Gromov-Witten invariants which arise.

\subsection{Ordered product formulas}
Let $\foD$ be a scattering diagram with $n$ lines
through the origin,
$$\foD = \{ (\fod_i,f_i)\ | \ 1\leq i \leq n \}\ .$$
Let $\Scatter(\foD)$ be the unique minimal scattering
diagram obtained by adding rays to
$\foD$ for which the path ordered product around the
origin is trivial.
Let $m'_\out\in M$ be a primitive vector. 
The function
$f_{m'_\out}$
is determined by the method used to prove
Theorems \ref{hhhp} and \ref{hhhpx}.

Let $m_1,\ldots, m_n$ be the primitives corresponding to the
lines of $\foD$.
Either $m'_\out$ is distinct from the $m_i$ or $m'_\out=m_k$
for some $k$.
Consider the
 $n=3$ case with functions  
$$f_1(s_\bullet^\bullet),\ 
f_2(t_\bullet^\bullet),\ f_3(u_\bullet^\bullet)$$
of the form \eqref{zzxb}.
If $m'_\out$ is distinct, then
\begin{equation}\label{gvgg}
\log f_{m'_\out} = \sum_{k=1}^\infty 
\sum_{\bf{G}=({\bf G}_1, {\bf G}_2, {\bf G}_3)}
        k\ N_{\bf{m}}[{\bf{G}}]\ s^{{\bf G}_1}\ t^{{\bf G}_2} 
\ u^{{\bf G}_3} \
z^{km'_\out}
\end{equation}
where the sum is over all graded partitions ${\bf G}_i$
of lengths $(\ell_i^1, \ldots, \ell_i^{d_i})$ satisfying
$$|{\bf G}_1|m_1 + |{\bf G}_2| m_2 
+ | {\bf G}_3|m_3
= k m'_\out.$$
The same result holds for all $n$.

For the degenerate case $m'_\out=m_k$, 
the definition of the invariant 
$N_{\bf{m}}[{\bf{G}}]$ must be changed slightly.
The only difference is the outgoing contact point $Q$ is
placed on the divisor $D_k^o$ instead of $D^o_\out$
(as discussed in Section \ref{degcas}). Then, equation
\eqref{gvgg}
holds as written. We leave the straightforward
details in the degenerate
case to reader.{\footnote{For the commutator formulas
of Theorem \ref{hhhp} and \ref{hhhpx}, nontrivial
degenerate cases do not appear. The full
arguments have been given there.}}

\subsection{Higher genus}
The higher genus analogues of the genus 0 invariants
$N_{\mathbf m}[{\mathbf G}]$ are not hard to construct.
Let 
$\overline{\foM}_g(X^o[{\bf G}]/D_\out^o)$ be the moduli space
of genus $g$ stable relative maps representing the
class $\beta_{\bf G}$ defined as in the genus 0 case. There are
now two difficulties: 
\begin{enumerate}
\item[(i)] $\overline{\foM}_g(X^o[{\bf G}]/D_\out^o)$ is not proper,

\item[(ii)] $\overline{\foM}_g(X^o[{\bf G}]/D_\out^o)$ 
is of virtual dimension $g$.
\end{enumerate}
The issues are resolved simultaneously by defining
$$N_{\mathbf m}^g[{\mathbf G}] =
\int_{[\overline{\foM}_g(X^o[{\bf G}]/D_\out^o)]^{vir}} (-1)^g \lambda_g$$
where $\lambda_g$ is the top Chern class of the 
Hodge bundle.{\footnote{See \cite{mp} for a parallel
definition of higher genus invariants in
the case of $K3$ surfaces.}} Limits out of the moduli space
$\overline{\foM}_g(X^o[{\bf G}]/D_\out^o)$ lead to loops in the
domain curve by the proof of Proposition \ref{pprr1}. However, the
class $\lambda_g$ vanishes on the locus of curves with loops.
Hence, $N_{\mathbf m}^g[{\mathbf G}]$ is well-defined.
A very interesting question is whether the 
relative invariants
$N_{\mathbf m}^g[{\mathbf G}]$
can be related to
the tropical vertex group.

\section{BPS state counts}
\label{bps}

\subsection{Log Calabi-Yau}
Let $S$ be a nonsingular surface and let $D\subseteq S$
be a nonsingular divisor.
The pair $(S,D)$ is {\em log Calabi-Yau} with respect to
$0\neq \beta\in H_2(S,{\mathbb{Z}})$ if
$$D\cdot \beta = c_1(S) \cdot \beta.$$
Two basic examples are:
\begin{enumerate}
\item[$\bullet$]
$(X^o[{\bf{P}}],D_\out^o)$,
constructed in Section \ref{sbu},
is log Calabi-Yau with respect to the class $\beta_{{\bf P}}$.
\item[$\bullet$]
$({\mathbb{P}}^2, E)$, where $E$ is a nonsingular
cubic, is log Calabi-Yau with respect to every class
$\beta\in H_2({\mathbb{P}}^2, {\mathbb{Z}})$.
\end{enumerate}

The moduli space $\overline{\foM}(S/D,w)$
 of genus $0$ stable relative
maps to $S/D$ of class $\beta$ and full tangency of order
$w=D\cdot \beta$ at a single unspecified point of $D$ is 
of virtual dimension 0. Let
$$N_S[w]\in {\mathbb Q}$$
be the associated relative Gromov-Witten invariant.
If 
\begin{equation}\label{jjjny}
\iota: P \rightarrow S
\end{equation}
is a rigid element of $\overline{\foM}(S/D,w)$, we can ask
what is the contribution of $d$-fold multiple covers of
$P$ to the Gromov-Witten invariant
$N_S[dw]$ in class $d\beta$?

\subsection{Multiple cover contributions}
We pursue here multiple cover calculations and
BPS state count definitions following the
perspective of \cite{pan1,pan}. In particular,
we assume the map $\iota$ of \eqref{jjjny} is as
well-behaved as possible.

Let $\infty \in P$ be the point of contact with $D$.
Let 
$$\overline{\foM}(P/\infty,d)^* \subseteq \overline{\foM}(S/D,dw)$$
be the locus of genus 0 stable relative maps 
$$(C,Q) \rightarrow P/\infty \rightarrow S/D$$
which factor as $d$-fold
covers of $P$. The moduli space $\overline{\foM}(P/\infty,d)^*$
is a nonsingular Deligne-Mumford stack of dimension $d-1$.
The superscript $*$ is used since the locus differs
slightly from the standard moduli space of stable relative maps
$\overline{\foM}(P/\infty,d)$. The reason is the $w$-tangency of $P$
with $D$ forces the ramification orders of maps over the destabilizations
to all be divisible by $w$. We leave the details here
for the reader.

The contribution $M_P[d]$ of $d$-fold multiple covers
of $P$ to $N_S[dw]$ is defined by
$$M_P[d] = \int_{[\overline{\foM}(P/\infty,d)^*]} e(B_d)$$
where $B_d$ is the obstruction bundle of rank $d-1$.
On the open locus of $\overline{\foM}(P/\infty,d)^*$
consisting of maps 
$$\varphi:(C,Q) \rightarrow P/\infty$$
with {\em no}  destabilizations of the
target, the obstruction space is
\begin{equation}\label{ghtt}
H^1\left(C, \varphi^*({\text{Norm}}_{S/P})(-(dw-d)Q)\right).
\end{equation}
By adjunction, the degree of the normal bundle
 ${\text{Norm}}_{S/P}$ is $w-2$.
Hence, the degree of $\varphi^*({\text{Norm}}_{S/P})(-(dw-d)Q)$ is
$-d$ and the obstruction space \eqref{ghtt} has rank $d-1$.
A description of the obstruction space for relative
maps can be found in \cite{gv}.

\begin{proposition} We have
$$M_P[d]= \frac{1}{d^2}\binom{d(w-1)-1}{d-1}\ .$$ 
\end{proposition}
\vspace{10pt}

If $w=1$, then $\binom{-1}{d-1}=(-1)^{d-1}$
by definition and the contribution
$$M_P[d]= \frac{(-1)^{d-1}}{d^2}\ ,$$
 specializes
to the genus 0 cap of \cite{Brypan}.

\vspace{10pt}

\proof
The $\CC^*$-action on $P$ fixing $\infty$ lifts to a
$\CC^*$-action on 
the moduli space $\overline{\foM}(P/\infty,d)^*$.
Once the lifting of $\CC^*$ to ${\text{Norm}}_{S/P}$
is chosen, a lifting of $\CC^*$ to the obstruction bundle
$B_d$ is determined
by the characterization of the obstruction space \cite{gv}.
Let $[1,-1]$ be the tangent weights of $\CC^*$ at the fixed points
$0,\infty \in P$.
We chose a lifting of 
$\CC^*$ to ${\text{Norm}}_{S/P}$ by specifying
 fiber weights 
$[w-2,0]$ over the respective fixed points.

We compute the integral $M_P[d]$ via the
 Bott residue formula.
There are many $\CC^*$-fixed loci in the moduli space
$\overline{\foM}(P/\infty,d)^*$.
 However,
the above linearization leads to the vanishing
 of {all}
contributions except for the single $\CC^*$-fixed Galois
cover
\begin{equation}\label{vvvhh1}
\varphi: (C,Q) \rightarrow P/\infty
\end{equation}
of degree $d$.
If $w=1$, the vanishing is the same as in the proof of
Theorem 5.1 of \cite{Brypan}.

If $w>1$, a different argument is needed.
Consider a $\CC^*$-fixed locus of 
$$\foL\subseteq \overline{\foM}(P/\infty,d)^*$$
 for which the target
is destabilized. 
For
$$[\varphi: (C,Q) \rightarrow P/\infty]\in \foL,$$
let $C',C''\subseteq C$ be the subcurves
mapped by $\varphi$ to the original and  destabilizing components of the
target respectively. 
The
$\CC^*$-action on the pull-back to $C'$
of ${\text{Norm}}_{S/P}$ is nontrivial. However,
 since $w>1$, 
$$H^1(C',\varphi^*({\text{Norm}}_{S/P}))=0.$$
The $\CC^*$-action on the pull-back to $C''$
of ${\text{Norm}}_{S/P}$ is trivial by our
choice of lifting.
The
 $\CC^*$-action on the destabilizing
components of the target is trivial.
The point $Q$ must map to a destabilizing
component of the target,
hence the $\CC^*$-action on ${\mathcal{O}}_C(-Q)$
is trivial.
By examining the obstruction space \cite{gv}, we conclude
the $\CC^*$-action on $B_d$ is {\em trivial} over $\foL$.
By dimension considerations, the contribution of $\foL$
 vanishes in the Bott residue formula ---
specifically the Euler class $e(B_d)$ is 0 when restricted to
$\foL$. 
The Galois cover \eqref{vvvhh1} is the unique
$\CC^*$-fixed locus for which the target is not
destabilized.

We compute the contribution of the Galois cover to the
Bott residue formula.
The weights of $\CC^*$ on the fibers of
$\varphi^*({\text{Norm}}_{S/P})(-(dw-d)Q)$
over the respective fixed points on $C$ are $[w-2,w-1]$.
The weights of $\CC^*$ on 
$H^1\left(C, \varphi^*({\text{Norm}}_{S/P})(-(dw-d)Q)\right)$
are
$$ \prod_{i=1}^{d-1} \frac{dw-d-i}{d} $$
The weights of $\CC^*$ on 
the tangent space to
$[\varphi]$ in 
$\overline{\foM}( {P}/\infty,d)^*$ are
$$\frac{\prod_{i=1}^{d}\frac{i}{d}}{\frac{1}{d}}$$
where the bottom factor is obtained from 
reparameterization over $0\in P$.
By the Bott residue formula,
$M_P[d]$ equals the ratio of the above weights, 
$$M_P[d] = \frac{1}{d} \frac{
\prod_{i=1}^{d-1} \frac{dw-d-i}{d} 
}
{
\frac{\prod_{i=1}^{d}\frac{i}{d}}{\frac{1}{d}}
}
= 
\frac{1}{d^2}\binom{d(w-1)-1}{d-1}\ ,
$$
together with a stack
automorphism factor of $\frac{1}{d}$. 
\qed

\subsection{Conjectures}
Let $(S,D)$ be a log Calabi-Yau pair with respect to
a primitive class $\beta\neq 0$.
Let $w=D\cdot \beta$ as before. 
Consider the generating series
$$N_{S} = \sum_{k=1}^\infty N_S[kw]\ q^k.$$
Using the multiple cover calculation, we can write
\begin{equation}\label{nt341}
N_{S} = \sum_{k=1}^\infty n_S[kw] \sum_{d=1}^\infty  
\frac{1}{d^2}\binom{d(kw-1)-1}{d-1} \ q^{dk}
\end{equation}
for unique numbers $n_S[kw]\in {\mathbb{Q}}$.
Equation \eqref{nt341} {\em defines} the $n_S[kw]$.

\begin{conjecture} The $n_S[kw]$ are integers for all $k\geq 1$.
\label{iinntt}
\end{conjecture}

Extracting integers from genus 0 Gromov-Witten theory
by removing multiple cover contributions is a basic
idea in the subject --- first pursued in the study of 
genus 0 curves on the quintic 3-fold. By 
the string theoretic
work of Gopakumar and Vafa \cite{GoVafa1,GoVafa2}, the resulting integers
can often be interpreted as BPS state counts in 
related theories. 
We interpret $n_S[kw]$ here as the associated BPS count.
Unfortunately, in almost every case, integrality
statements of the form of Conjecture \ref{iinntt}
are not provable by existing techniques.{\footnote{One
exception is the Fano 3-fold case settled in
\cite{Z}. Calabi-Yau cases such as ours here
are more difficult.}}

Conjecture \ref{iinntt} applies to the geometries
$(X^o[{\bf P}], D^o_\out)$ associated to the functions
$f_{{m}'_\out}$ in 
Theorem \ref{hhhp}.  
Kontsevich and Soibelman conjectured \cite{ks2} an equivalent
integrality{\footnote{This conjecture, in
our language, would posit a multiple cover
contribution of the form $(-1)^{d-1}/d^2$ independent
of $w$. As we have seen, the multiple
cover contributions are more subtle (and depend upon $w$).
But, the associated integrality is the same.}}
 for the functions 
$f_{m'_\out}$ associated to such commutators, and 
a proof via quiver techniques has been recently provided in
\cite{Rein}.
Conjecture \ref{iinntt} applies to other quite
different situations as well. The most interesting
case is perhaps the log Calabi-Yau pair $(\mathbb{P}^2, E)$
studied in \cite{Ga,Ta}. In all the examples
we have considered, there is good numerical evidence
supporting Conjecture \ref{iinntt}.

Conjecture \ref{iinntt} does {\em not} apply as stated to
the orbifold geometries $(X^o[{\bf G}]/D^o_\out)$
of Theorem \ref{hhhpx}. The orbifold structure
leads to more complicated multiple cover contributions
which we have not yet calculated.

\subsection{Examples}

Consider the  $\ell_1=\ell_2=3$ case of the 
commutator \eqref{jjtr} of the Introduction. The
corresponding scattering diagram is discussed
in Example \ref{basicexample} of Section \ref{pop}. Focus
on the function attached to the line of slope $1$. 
By direct calculation in the tropical vertex group, we find
\begin{multline*}
\log f_{\out}=
9(t_1t_2xy)+
2\cdot {63\over 4}(t_1t_2xy)^2+
3\cdot 55(t_1t_2xy)^3  \\
+
4\cdot {4095\over 16}(t_1t_2xy)^4+
5\cdot {100947\over 25}(t_1t_2xy)^5+\cdots
\end{multline*}

Consider $\PP^2$ with
the coordinate axes $D_1,D_2$ and $D_{\out}$ and
pick points 
$$x_{11},x_{12},x_{13} \in D_1, \ \
 x_{21}, x_{22}, x_{23}\in D_2 \ .$$
Sums of the relative Gromov-Witten invariants 
$N_{\bf m}[({\bf P}_1, {\bf P}_2)]$ are
determined by Theorem 5.4. For example, 
\[
\sum_{|{\bf P}_1|=1,\ |{\bf P}_2|=1}
N_{\bf m}[({\bf P}_1, {\bf P}_2)]=9.
\] 
We can easily interpret the answer
in the following way.  Given a choice of
one of the three points on $D_1$ and one of the three points on $D_2$,
there is precisely one line through these two points, which of course
is maximally tangent to $D_{\out}$. There are nine such configurations, 
hence
the correct answer is $9$.

The next coefficient of $f_\out$ yields
\[
\sum_{|{\bf P}_1|=2,\ |{\bf P}_2|=2}
N_{\bf m}[({\bf P}_1,{\bf P}_2)]={63\over 4}.
\] 
The double covers of the lines
mentioned above count for $-9/4$, leaving a contribution of $63/4+9/4=18$
from conics passing through two of the three points on $D_1$ and two of
the three points on $D_2$. Indeed, given any choice of two points each
on $D_1$ and $D_2$, there are two conics through these four points
tangent to $D_{\out}$.

The third coefficient  is more interesting. We have
\[
\sum_{|{\bf P}_1|=3,\ |{\bf P}_2|=3}
N_{\bf m}[({\bf P}_1,{\bf P}_2)]=55.
\] 
The contribution from triple covers of the lines is 
$9/3^2=1$, hence we expect 54 non-multiply covered cubics.
A refined scattering diagram calculation reveals more specifically
the following numbers:
\begin{eqnarray*}
N_{\bf m}[1+1+1,1+1+1]&=&18\\
N_{\bf m}[2+1+0,1+1+1]&=&3
\end{eqnarray*}
The first number can be interpreted as the number of nodal cubics
passing through all six of the chosen points and maximally tangent
to $D_{\out}$, while the second number can be interpreted as
the number of nodal cubics passing through $x_{11}$ twice (so that
the node is at $x_{11}$), passing through $x_{12}, x_{21}, x_{22}$ and $x_{23}$
once, and again maximally tangent to $D_{\out}$. There are a total
of $12$ partitions ${\bf P}=({\bf P}_1,{\bf P}_2)$ 
involving the numbers $2,1,0,1,1,1$ of
this sort, so the total accounting is
\[
54=18+12\times 3.
\]
In the above cases, the Gromov-Witten invariants (corrected
for multiple covers)
solve straight counting problems.

\section*{Appendix: Tropical/holomorphic counts}

The purpose of this Appendix is to discuss the modifications of
\cite{nisi} necessary to obtain a proof of Theorem~\ref{N=Ntrop}. 
We do not try to be self-contained and rather just indicate what
has to be changed.

First, some conventions of \cite{nisi} conflict with the notation
adopted in the present paper. Most importantly, the roles of $M$ and
$N$ are reversed. The reason is that fans and tropical curves
traditionally live in $N$, but the tropical vertex group naturally
acts on polyomial rings with exponents on the dual lattice $M$.
Moreover, $n$ denotes the rank of $N$ in \cite{nisi}, while now we
work in dimension~$2$ and $n$ denotes the number of incoming
directions. In \cite{nisi}, our toric variety $X$ is written
$X(\Sigma)$ and $X$ denotes the total space of a toric degeneration.
Another irrelevant difference is that in \cite{nisi} we do tropical
geometry over $\QQ$ while here we work over $\RR$. In the Appendix, we
follow the notation of \cite{nisi} except for swapping $N$ and $M$.

The degree of the tropical curves to be considered is fixed by the
number and directions of the incomming edges. Rather than
imposing incidence of a marked edge with an affine subspace as in
 Definition~1.3 of \cite{nisi}, we constrain the incoming unbounded edge
$E_{ij}$ by the choice of an element $m_{ij}\in M_\QQ/\QQ m_i$. The latter
is equivalent to the condition 
$$h(E_{ij})\subseteq \fod_{ij}=
m_{ij}+\QQ m_i$$ of
 Definition~\ref{def.Ntrop}. 
The tuple $\mathbf A=(m_{ij})$ with $m_{ij}\in M_{\QQ}/\QQ m_i$
determines a {\em constraint}.
A marked tropical curve
$(\Gamma,(E_{ij}),h)$ \emph{matches} the constraint $\mathbf{A}$ if
for all $i,j$
\[
h(E_{ij})=\{ m_{ij}\}\quad\text{in }M_\QQ/\QQ m_i.
\]
The finiteness and transversality results in Section~2 of \cite{nisi}
carry over without difficulty, but are partly already contained in
Theorem~\ref{onetoonecorrespondence}. The
{\em gluing map} (4) in \cite{nisi} now reads   
\begin{equation}\label{transversality sequence}
\begin{split}
\quad\quad\Phi: \operatorname{Map}(\Gamma^{[0]},M_\QQ)\ &\lra\ 
\prod_{E\in\Gamma^{[1]} \setminus \Gamma_\infty^{[1]}} M_\QQ/\QQ
u_{(\partial^-E,E)} \times \prod_{i,j} M_\QQ/ \QQ m_i,\\
h\ &\longmapsto\ \big((h(\partial^+E)-h(\partial^-E))_E,
h(\partial^-E_{ij})-m_{ij}\big).
\end{split}
\end{equation}
\smallskip

In Section~3 of \cite{nisi} we need to adapt the treatment of
incidence points, which in the present situation lie on the toric
boundary. Let
$$\pi:X=X(\widetilde\Sigma_\P)\to\AA^1$$
be the toric degeneration defined by an integral polyhedral
decomposition $\P$ of $M_\QQ$, and let $X_t=\pi^{-1}(t)$. Let
$D\subseteq X$ be the union of those toric divisors of $X$ not
contained in $X_0$, that is, corresponding to rays of
$\widetilde\Sigma_\P$ contained in $M_\QQ\times\{0\}$.  The
intersection $D_t=D\cap X_t$ for $t\neq 0$ is the toric boundary of
$X_t$, while $D_0\subseteq X_0$ is the union of those 1-codimensional
toric strata not contained in the singular locus of $X_0$. We thus
consider $D_0$ as the toric boundary of $X_0$.

The asymptotic fan $\Sigma_\P$ of $\P$ is the fan defining $X_t$ for
any $t\neq0$. Thus if $\omega =\QQ_{\ge0}\cdot  u$ is a ray of
$\Sigma_\P$, then $\omega\times\{0\}$ is a ray of $\widetilde\Sigma_\P$.
The fan describing the associated toric divisor 
$\widetilde D_\omega\subseteq D$ consists of the images under the projection
$$M_\QQ\times\QQ\to (M_\QQ/\QQ u)\times\QQ$$ 
of the cones $\sigma\in \widetilde\Sigma_\P$ containing
$\omega\times\{0\}$.  This is again a fan of cones over the cells of a
polyhedral decomposition $\P_\omega$, now of $M_\QQ/\QQ u$. The
associated toric degeneration is $\pi|_{\tilde D_\omega}$. The
vertices of $\P_\omega$, or equivalently the irreducible components of
$$D_\omega:= \widetilde D_\omega \cap X_0,$$  are in bijective
correspondence with the unbounded edges of $\P$ in direction $\omega$.
If $\widetilde\omega$ is such an unbounded edge and $D_{\widetilde
\omega}\subseteq D_\omega$  the corresponding irreducible component,
the linear space generated by $\widetilde \omega\times\{1\}$ in
$(M_\QQ/\QQ u)\times\QQ$ defines a 1-dimensional subtorus
$\GG_{\widetilde \omega}$ of $\GG((M/\ZZ u)\times\ZZ)$, the stabilizer
of $D_{\widetilde \omega} \subseteq \widetilde D_\omega$.  By
Corollary~3.8 of \cite{nisi}, given a closed point $Q\in D_t$ for
$t\neq 0$ and assuming $\P$ integral, the closure of the orbit
$\GG_{\tilde \omega}\cdot Q$ is a section $\widetilde Q$ of 
$$\widetilde D_\omega\to \AA^1$$ with $\widetilde Q\cap
X_0\subseteq D_{\widetilde\omega}$. 
Summarizing, the choice of $l_i$
unbounded edges in direction $u=-m_i$ in a polyhedral decomposition
with asymptotic fan $\Sigma$ readily defines a degeneration of our
incidence points $Q_{i1}, \ldots, Q_{i l_i}$ to points $Q^0_{ij}\in
X_0$ on disjoint toric strata of $\widetilde D_{\omega_i}$, $\omega_i=
-\QQ_{\ge 0} m_i$, the degeneration of the toric divisor $D_i\subseteq
X_t= X(\Sigma)$.
\smallskip

For a general constraint $\mathbf A$, the affine map
\eqref{transversality sequence} is an isomorphism. If $\foD\in\NN$ is
the index of the corresponding inclusion of lattices, Proposition~5.7
of \cite{nisi}
about the existence of exactly $\foD$ isomorphism classes of maximally
degenerate curves matching the incidence conditions on $X_0$
works as before. The only difference is that the $ij$-th incidence
condition is now a torsor under $\GG(M/ \ZZ m_i)$. 
\smallskip

In the deformation theory of Section~7 of \cite{nisi},
 the discussion of
the situation at the toric boundary is somewhat hidden. In the
notation of \cite{nisi}, the above discussion provides the degeneration
of toric boundaries $\underline D\subseteq\underline X$. Now there are
functions on $X$ vanishing along $\underline D$, but not on
$\underline{X_0}$. These force the introduction of special (marked)
points in the log-structure of $C_0$, and the order of tangency of
$\underline{\varphi_0}$ with $\underline D$ at such a point $x\in
\underline{C_0}$ is fixed by the induced map of {\em ghost sheaves}
\[
\underline{\varphi_0}^*\overline{\M}_{X}\lra \overline{\M}_{C_0}.
\]
In fact, $\overline{\M}_{X,\underline{\varphi_0}(x)}$ has a direct
summand $\NN$ generated by a local equation for $\underline D$, and
similarly $\overline{\M}_{C_0,x}$ has a direct summand $\NN$ coming
from a local equation for the marked point. The map on these direct
summands is multiplication with the order of tangency. Log deformation
theory preserves this map and hence considers deformations with fixed
tangency conditions from the outset. Thus there is nothing to be
changed here.

The only difference is again the discussion of the incidence
conditions via the transversality argument in 
Proposition~7.3 of \cite{nisi}.  
In our situation, the point is to show surjectivity
of the evaluation map
\begin{eqnarray}\label{transversality map}
H^0(\shN_{\varphi_0})\lra
\prod_{i,j} T_{X/{\AA^1},\underline{\varphi_0}(x_{ij})} \big/
D\varphi_0(T_{C_0/O_0,x_{ij}}),
\end{eqnarray}
just as in (9) of \cite{nisi}. This is the same as the second
component
\[
\operatorname{Map}(\Gamma^{[0]},M)\ \lra\ 
\prod_{i,j} M/ \QQ m_i
\]
of \eqref{transversality sequence}, tensored with $\CC$, so is
surjective for tropical curves that are general in the sense of
Definition~2.3 of \cite{nisi}.

Therefore, fixing $\mathbf A$ general, the arguments of Section~8 of
\cite{nisi} produce a bijective
 correspondence between holomorphic
curves contributing to $N_{\mathbf m}^\hol (\mathbf w)$ in $X_t$ for
small $t$ and certain stable log-maps to $X_0$, for a fixed
degeneration $$X\to\AA^1.$$
 Each stable log map yields a tropical curve
contributing to $N_{\mathbf m}^\trop(\mathbf w)$, computed with the
asymptotics provided by $\mathbf A$. Conversely, a tropical curve
$(\Gamma,\mathbf E,h)$ has $\foD(\Gamma,\mathbf E,h)\cdot w(\Gamma)$
stable log maps associated to it. Here $\foD(\Gamma,\mathbf E,h)$ is
the lattice index associated to \eqref{transversality sequence}, and
$w(\Gamma)$ is the product of all weights of bounded edges. Thus
$(\Gamma,\mathbf E,h)$ contributes
\[
w(\Gamma)\cdot \foD(\Gamma,\mathbf E,h)
\]
to $N_{\mathbf m}^\hol(\mathbf w)$. Now $\foD(\Gamma,\mathbf E,h)$ is
the same as the lattice index $\foD(\Gamma,\mathbf E,h,\mathbf P)$  of
\cite{nisi} for imposing pointwise incidence conditions $\mathbf
P=(P_{ij})$ on $h(E_{ij})$. Proposition~8.8 of \cite{nisi} thus
implies
\[
w(\Gamma,\mathbf E)\cdot \foD(\Gamma,\mathbf E,h)
=\prod_{V\in\Gamma^{[0]}} \Mult_V(h),
\]
where we used the notation of Definition~\ref{multiplicitydefinition}
and $w(\Gamma,\mathbf E)= w(\Gamma)\cdot\prod_{i,j} w_{ij}$ in the
present case.\footnote{The statement of Proposition~8.8 of \cite{nisi}
with $w(\Gamma)$ rather than $w(\Gamma,\mathbf E)$ is wrong, and in
fact, $w(\Gamma,\mathbf E)$ is also needed there for the claimed equivalence
with the mutitplicity of \cite{Mk} . The problem is an incorrect
verification of the base case of the induction. The rest of the proof
remains the same.} Hence,
\[
N_{\mathbf m}^\hol (\mathbf w)= \frac{1}{\prod_{i,j}w_{ij}}
\sum_{\{(\Gamma,\mathbf E,h)\}} \prod_{V\in\Gamma^{[0]}} \Mult_V(h)
=
\frac{
N_{\mathbf m}^\trop(\mathbf w)}
{\prod_{i,j}w_{ij}}
,
\]
finishing the proof of Theorem~\ref{N=Ntrop}.
\qed


\end{document}